\documentclass[12pt]{amsart} 
\usepackage{amsmath,amssymb,amscd,array}

\usepackage{mathptmx}
\usepackage{graphics}

\def\url#1{\expandafter\string\csname #1\endcsname}
\usepackage{url}

\usepackage{tabularx,mdwtab}

\usepackage{bbm}

\newcommand {\foo}{{\mathfrak{oo}}}
\def\mmat #1,#2,#3,#4,{\text{\small\arraycolsep=3pt $
\begin{pmatrix}#1&#2\\#3&#4\end{pmatrix}$}}

\newcommand {\pd}[1] {\partial_{#1}}

\usepackage{DLdef1}

\newcommand {\fbr}{{\mathfrak{br}}}

\newcommand {\fwk}{{\mathfrak{wk}}}
\newcommand{\un}{\underline{N}}
\newcommand{\del}{\partial}
\rmnameii{etr}{etr}
\rmnameii{evv}{ev}

\begin{document}

\title[Classification of simple Lie superalgebras]{Classification of simple Lie superalgebras in characteristic $2$}

\author{Sofiane Bouarroudj${}^a$, Alexei Lebedev${}^b$,\\
Dimitry Leites${}^{a,c}$, Irina
Shchepochkina${}^d$}

\address{${}^a$New York University Abu Dhabi,
Division of Science and Mathematics, P.O. Box 129188, United Arab
Emirates; sofiane.bouarroudj@nyu.edu \\ ${}^b$Equa
Simulation AB, R{\aa}sundav\"agen 100, Solna, Sweden; alexeylalexeyl@mail.ru\\
${}^{c}$Department of Mathematics, Stockholm University, SE-106 91
Stockholm, Sweden; dl146@nyu.edu;
mleites@math.su.se\\
${}^d$Independent University of Moscow, Bolshoj Vlasievsky per, dom
11, RU-119 002 Mos\-cow, Russia; irina@mccme.ru}



\begin{abstract} All results concern characteristic 2. Two procedures that to every simple Lie algebra
assign simple Lie superalgebras, most of the latter new, are offered. We
prove that every simple finite-dimensional Lie superalgebra is
obtained as the result of one of these procedures, so we classified
all simple finite-dimensional Lie superalgebras modulo non-existing
at the moment classification of simple finite-dimensional Lie
algebras. 

This result concerns Lie
superalgebras considered naively, as vector spaces.
To obtain classification of simple Lie superalgebras in the category
of supervarieties, one should list
the non-isomorphic deforms (results of deformations) with odd parameter. This problem is open
bar several examples described in arXiv~0807.3054.

For Lie algebras,
 in addition to the known ---
``classical" --- restrictedness, we introduce a (2,4)-structure on the two non-alternating series: orthogonal  and
of Hamiltonian vector fields. For Lie superalgebras, the classical
restrictedness of Lie algebras has two analogs: $(2|4)$- and
$(2|2)$-structures; one more analog --- a $(2,4)|4$-structure on Lie superalgebras is
the analog of  (2,4)-structure on Lie algebras.
\end{abstract}

\thanks{S.B. was partly supported by the grant AD 065 NYUAD. D.L. is thankful to l'IHES for several summer months spent in it in 1990s. The results were delivered, in particular, at the conference in honor of A.~Kirillov, Reims, 2017.
We are thankful to S.~Skryabin,
P.~Zusmanovich and A.~Krutov for help.}


\markboth{\itshape Sofiane Bouarroudj\textup{,} Alexei
Lebedev\textup{,} Dimitry Leites\textup{,} Irina
Shchepochkina}{{\itshape
Classification of simple Lie superalgebras in characteristic $2$}}

\maketitle

\thispagestyle{empty}

\section{Introduction}\label{Intro}

The Lie algebras and superalgebras we consider in this paper are
finite-dimensional (except for \S\ref{Z2gra}) over an algebraically
closed field $\Kee$ of characteristic $p>0$, mostly, for $p=2$. (Observe that neither finite-dimensionality of algebras,
nor algebraic closedness of the ground field is needed in most of
the constructions below.) In examples, we need (generalized) Cartan prolongation (for more examples, see \cite{BGLLS2}) and Lie (super)algebra with indecomposable Cartan matrix over $\Kee$ classified in \cite{BGL2}.

\ssec{What restricted Lie algebra in characteristic 2 is} In 2005, P.~Deligne wrote several
comments to a paper by the two of us, see his Appendix in \cite{LL}.
In particular, a part of his advice is (in our words):
``Over $\Kee$, to classify ALL simple Lie (super)algebras and their
representations are, perhaps, not very reasonable problems, and
definitely very tough; investigate first the \textbf{restricted}
case: it is related to geometry, meaningful and of interest". 

But what is restrictedness if $p=2$? For possible answers, see \S\ref{Srestric}, where a reason for
difference in answers for $p=2$ and $p>2$ is conjecturally attributed to the fact~\eqref{reason}.

Having read \cite{LL}, several algebraists argued that in certain problems of interest to experts in the field, non-restricted Lie algebras are important, see \cite{Kos}. Non-restricted Lie algebras are also needed even to describe the simple restricted algebras, see \cite{BW, SF, S, BGP, Sk}. So in our quest for simple Lie (super)algebras we do not restrict ourselves to restricted algebras (whatever they are, cf. \cite{WZ}).

\ssec{On classification of simple Lie algebras over $\Kee$} \underline{For $p>7$}, Kostrikin and Shafarevich suggested (ca 1966) a construction of simple Lie algebras and a conjectural list of all restricted simple Lie algebras thus obtained,  see \cite{KS}. This list lucidly describes $\Zee$-graded Lie algebras, but their deforms were clearly described and classified only much later, e.g., see \cite{Sk} and \cite{W}.

\underline{For $p>7$}, Block and Wilson proved the restricted version of the KSh-conjecture, \cite{BW}. This classification is implicit when dealing with deforms: simple deforms of the divergence-free algebras (\cite{W}) and of hamiltonian type algebras  (\cite{Sk}) were classified only several years after  \cite{BW} was published;  explicit formulas of the $p$-structure were obtained only recently 
in the divergence-free case and (despite huge computational difficulties) in some of Hamiltonian cases, see \cite{BKLS}.

\underline{For $p>3$},  the KSh-conjecture amended by adding Melikyan algebras (for $p=5$) was
proved only recently, see \cite{S,BGP}. Classification of filtered deforms of Hamiltonian Lie algebras due to Skryabin, see \cite{Sk}, is a very important ingredient difficult to prove.

After his teacher (A.I.Kostrikin)
died, A.~Dzhu\-ma\-dil\-daev wrote the paper \cite{KD}, where he
suggested to improve the Kostrikin-Shafarevich construction by
identifying certain Lie algebras as ``standard" examples, every simple
Lie algebra being either on the list of ``standard" algebras, or a \textit{deform} (i.e.,
the result of a deformation) of one of the standard objects.

\underline{For $p\leq 5$}, \textbf{the stock of ``standard"
examples should contain not only simple Lie alge\-bras}, but al\-so their nontrivial
central extensions (or even semi-simple Lie algebras, see \cite{SkT1}):
the (family of) Melikyan algebras are described in \cite{KD}
as deforms of certain Poisson Lie algebras (in \cite{KD},
the proof of the fact that the Melikyan algebra is a deform of the
Poisson algebra is correct, although the Poisson algebra is called
Hamiltonian; so the computer verification, see \cite{MeZu}, is not
needed). The KSh-construction thus improved embraces $p=5$. 

In
\cite{KD}, Dzhumadildaev also claimed that the Ermolaev and
(unspecified) Skryabin algebras are deforms of examples obtained by
the KSh-construction. For Ermolaev algebras, this claim was
announced almost a decade earlier in \cite{KuJa}; for the full set of
explicit cocycles corroborating the claim of \cite{KuJa}, see
\cite{BGL4}. The importance of
investigating deforms, and isomorphism classes and central
extensions thereof, is now manifest; for classification of
deforms of simple Lie
algebras with symmetric root system (containing together with every root its opposite of the same
multiplicity) and of
rank $\leq 8$ (the vital cases in the inductive proof of the future classification), see \cite{BGL3}.

\underline{True deforms}.  Further study revealed the following serious obstacle seldom discussed before. In \cite{DK}, one of the
first papers on deformations of simple vectorial Lie algebras, it
was observed that although $H^2(\fg;\fg)\neq 0$, and the cocycles
representing the nontrivial classes of $H^2(\fg;\fg)$ are
integrable, the deforms corresponding to some or all of
these cocycles can be isomorphic to $\fg$. We call such cocycles,
and deformations corresponding to them, \textit{semitrivial}, see
\cite{BLW}; a wide class of them is characterized in \cite{BGLLS, BGL3}; for
examples for $p=0$, see \cite{Ri}. Thus, we are only interested in \textbf{true deforms}, not trivial or semi-trivial deforms, of \lq\lq standard\rq\rq\ examples.

\underline{For $p=3$}, in addition to the simple Lie algebras
obtained by the KSh-construction, and deforms thereof, there are
known other examples: Frank algebra, Ermolaev algebra, and several
rather mysterious Skryabin algebras; each
of these algebras is, actually, a member of a  family depending on the shearing
vector $\un$. These Lie algebras were interpreted (demystified) in
\cite{GL3}, where the number of parameters on which the vector $\un$
depends was corrected as compared with previous descriptions
\cite{S}. These examples from \cite{GL3} conjecturally complete the
stock of ``standard" objects for $p=3$.

\underline{For $p=2$},
conjecturally, the set of ``standard"
Lie algebras is the union of examples from
\cite{Ei, SkT1, BGL2,LeP, ILL, BGLLS,BGLLS1,BGLLS2}, and those that might be
obtained by methods of \cite{LSh} from nontrivial central extensions
and algebras of derivations listed in \cite{BGLL1}, and from still to be listed algebras with non-symmetric root system.

\sssec{On classification of simple Lie superalgebras} Lie superalgebras appeared not in high energy physics in 1974 or in solid state physics in 1980s, as we sometimes hear and read, but in 1930s in  topology. They appeared in the form of --- in modern terms --- super Lie rings, i.e., over $\Zee$, or over finite fields. For example, the homotopy groups constitute a super Lie ring with respect to the Whitehead multiplication. These examples are solvable, hence without nice structure like that of simple ones, and this fact delayed the study of modular Lie superalgebras.

\underline{The case of $p=0$} illustrates
how much more difficult classification of simple Lie superalgebras  is as compared
with that of simple Lie algebras of the same type.

The  finite-dimensional  simple Lie superalgebras were classified, thanks to several groups of
researchers, inside two years or so, see \cite{Kapp}. For a review not only
of classification (with due references to results of Kaplansky, Djokovich--Hochschield, and Scheunert--Nahm--Ritten\-berg), but of
impressive at that time results on the basics of representation theory, see \cite{K0}.

For classification of simple $\Zee$-graded Lie superalgebras of polynomial
growth a big chunk of the problem (in particular, of vectorial Lie superalgebras
with polynomial and Laurent coefficients) is solved on hundreds of pages, see
\cite{LSh,Sh5,Sh14, GLS1, HS} and \cite{K,K10, CCK,CKa} and references therein,
but the conjecture formulated in \cite{LSS} embracing all types of Lie superalgebras ($\Zee$-graded of polynomial
growth) is still open.

\underline{For $p>0$}, the classification of simple Lie superalgebras is even more
difficult than classification of simple Lie algebras for  $p>0$, the smaller $p$, the more difficult. And for any $p>0$ it is more difficult than classification of simple Lie superalgebras for $p=0$.
For a conjectural classification of simple Lie superalgebras, rather complicated even for $p>2$, not proved for any $p$, even for
restricted Lie superalgebras, see \cite{BGLLS1}. (The absence of classification is no wonder: even the classification of restricted simple Lie algebras for $p>3$ is implicit to this day:  \lq\lq The problem of restrictedness is approached. ...[But] the family of Hamiltonian algebras ... is not yet handable", see \cite[v.1, p.357]{S}.) 

So the main result
of this paper is unexpected: it reduces the list of ``standard"
examples needed to construct simple Lie superalgebras for $p=2$ ---
presumably the most difficult case of all classifications spoken above --- to the list of ``standard"
examples for simple Lie algebras.

\ssec{Summary of our results} In this note we briefly overview the classification problem of simple Lie (super)algebras over $\Kee$ and
formulate our results expounding \cite{LeD}, where A.~Lebedev briefly described the
following three phenomena existing for $p=2$ only:

1) The two methods, see subsec.~\ref{Method1} (queerification) and \ref{Method2} (\lq\lq method 2"),
producing simple Lie superalgebras from every simple Lie algebra.
Here we prove that every simple Lie superalgebra is obtained by one of these procedures.

2) On Lie algebras: in addition to the known, ``classical", restrictedness there
exists a $(2,4)$-structure, e.g., on the orthogonal Lie algebras $\fo(2n+1)$ and the
Lie algebras of non-alternating Hamiltonian vector fields $\fh_I(2n+1;\underline{\One})$, where $\underline{\One}:=(1,\dots,1)$ is the shearing vector of heights of indeterminates.

3) On Lie superalgebras: there are three analogs of the classical
restrictedness, namely, (a) the classical (direct) one, i.e., $(2|4)$-structure, (b) the
$(2|2)$-structure, (c) the
$(2,4)|4$-structure. Hereafter by \textit{restricted Lie
(super)algebra} we mean a classically restricted one.

In \S\ref{Srestric}, we consider various restrictednesses, with
more details and examples than in \cite{BGL2}.

In \S\ref{Sproof}, we describe the two
methods (queerification and ``method 2") producing
new simple Lie superalgebras from every simple Lie algebra if $p=2$.

In \S\ref{Snew}, we prove our main Theorem: \textbf{if $p=2$,
every simple finite-dimensional Lie superalgebra is obtained from a
simple Lie algebra by means of either queerification, or ``method
2"}. So we have obtained classification of simple Lie superalgebras modulo
classification of simple Lie algebras. Here ``Lie
superalgebra" is understood naively, as a $\Zee/2$-graded algebra
satisfying certain identities. \textit{To classify Lie
superalgebras considered in the category of supervarieties we need
to describe odd parameters of deformations of naively understood Lie
superalgebras}, see \cite{BGL3}. 

Observe a ``braking news" result: the list of $\Zee/2$-gradings of simple Lie
algebras --- a vital ingredient for ``method
2" --- is much longer for $p=2$ than for the namesakes  of these algebras (``the same", in a sense, algebras) considered for $p\neq 2$, see \cite{KrLe}.

In \S\ref{Sq} and \S\ref{SrelCartAndQ}, examples of queerifications
are given with details further clarifying the mechanism which does
not exist for $p\neq 2$.

In \S\ref{Z2gra}, examples illustrating numerous new simple Lie superalgebras obtained by
means of ``method 2" from vectorial Lie algebras are explicitly described.


\section{Several versions of restrictedness in characteristic 2}\label{Srestric}

Here is an expounded with new results version of the respective
sections from \cite{LeD, BGL2}.

\ssec{Restrictedness on Lie algebras}\label{ss-p-str} Let the ground field $\Kee$
be of characteristic $p>0$, and $\fg$ a~Lie algebra. For every $x\in
\fg$, the operator $(\ad_x)^{p}$ is a derivation of
$\fg$. If this derivation is an inner one, i.e., there is a map
(called \textit{$p$-structure}) ${}[p]:\fg\tto\fg, \ x\mapsto x^{[p]}$ such that
\begin{equation} \label{restricted-3}
[x^{[p]}, y]=(\ad_x)^{p}(y)\quad \text{~for any~}x,y\in\fg,
\end{equation} 
\begin{equation}\label{restricted-1}
(ax)^{[p]}=a^px^{[p]}\quad \text{~for any~}a\in\Kee,~x\in\fg,
\end{equation}
\begin{equation} \label{restricted-2}
(x+y)^{[p]}=x^{[p]}+y^{[p]}+\mathop{\sum}\limits_{1\leq i\leq
p-1}s_i(x, y) \quad\text{~for any~}x,y\in\fg,
\end{equation}
where $is_i(x, y)$ is the coefficient of $\lambda^{i-1}$ in
$(\ad_{\lambda x+y})^{p-1}(x)$,
then the Lie algebra
$\fg$ is said to be \emph{restricted} or \emph{having a
$p$-structure}.

\sssbegin{Remarks}\label{Uniq} 1) If the Lie algebra $\fg$ is centerless,
then the condition \eqref{restricted-3} implies  \eqref{restricted-1} and  \eqref{restricted-2}.

A $p$-structure on a given Lie algebra $\fg$ does not have to be unique; all $p$-structures on $\fg$ agree modulo center. Hence, on any simple Lie algebra, there is
not more than one $p$-structure.

2) According to \cite[Th. 2.3, p. 71]{SF}, the following condition, due to Jacobson, is sufficient for a Lie algebra $\fg$ to
have a $p$-structure: for a basis $\{g_i\}_{i\in I}$ of $\fg$,
there exist elements $g_i^{[p]}$ such that
\[ [g_i^{[p]}, y]=(\ad_{g_i})^{p}(y)\quad \text{~for any~}y\in\fg.
\]
\end{Remarks}

\sssec{Restricted modules}\label{resModP} A $\fg$-module $M$ over a
restricted Lie algebra $\fg$, and the representation $\rho$ defining
$M$, are said to be \textit{restricted} or having a
\textit{$p$-structure} if
\begin{equation*} \label{restri-3}
\rho(x^{[p]})=(\rho(x))^{p}\quad \text{~for any~}x \in\fg.
\end{equation*}

\ssec{Lie
superalgebras}\label{SSlieSuper} Naively, the definition of \textit{Lie superalgebra}
is the same for any $p\neq 2$. Let us point at the subtleties for $p=2$. For any $p$, a \textit{Lie superalgebra} is a
superspace $\fg=\fg_\ev\oplus\fg_\od$ such that the even part
$\fg_\ev$ is a Lie algebra, the odd part $\fg_\od$ is a
$\fg_\ev$-module (made into the two-sided one by
\textit{anti}-symmetry, i.e., $[y,x]=-[x,y]$ for any $x\in
\fg_\ev$ and $y\in
\fg_\od$), and  a \textit{squaring}
defined on $\fg_\od$ as a map $S^2(\fg_\od)\tto\fg_\ev$:
\begin{equation*}\label{squaring}
\begin{array}{c}
x\mapsto x^2\in\fg_\ev\quad \text{such that $(ax)^2=a^2x^2$ for any $x\in
\fg_\od$ and $a\in \Kee$, and}\\
{}[x,y]:=(x+y)^2-x^2-y^2\text{~is a bilinear form on $\fg_\od$ with values
in $\fg_\ev$.}
\end{array}
\end{equation*}
(This extra requirement on squaring is needed, say, over $\Zee/2$ where
not any quadratic form that vanishes at the origin yields a bilinear form $[\cdot,\cdot]$.)

The Jacobi identity involving odd elements takes the following form:
\begin{equation*}\label{JI}
\begin{array}{l}
~[x^2,y]=[x,[x,y]]\text{~for any~} x\in\fg_\od, y\in\fg_\ev,\\
~[x^2,x]=0\;\text{ for any $x\in\fg_\od$.}
\end{array}
\end{equation*}

For any Lie \textbf{super}algebra $\fg$, its \textit{derived algebras} are
defined to be (for $i\geq 0$)
\[
\fg^{(0)}: =\fg, \quad
\fg^{(i+1)}=\begin{cases}[\fg^{(i)},\fg^{(i)}]&\text{for $p\neq
2$,}\\
[\fg^{(i)},\fg^{(i)}]+\Span\{g^2\mid g\in\fg^{(i)}_\od\}&\text{for
$p=2$}.\end{cases}
\]

\ssec{The $p|2p$-structure or restricted Lie
superalgebra}\label{SSp2pStr}

For a Lie superalgebra $\fg$ of characteristic $p>0$, let the Lie
algebra $\fg_\ev$ be restricted and
\begin{equation} \label{restr3}
[x^{[p]}, y]=(\ad_x)^{p}(y)\quad \text{~for any~}x\in\fg_\ev,~y\in\fg.
\end{equation}
This gives rise to the map (recall that the bracket of odd elements
is the polarization of the squaring $x\mapsto x^2$)
\begin{equation*}\label{2p}
{}[2p]:\fg_\od\to\fg_\ev, ~~~ x\mapsto(x^2)^{[p]},
\end{equation*}
satisfying the condition
\begin{equation*}\label{2p2}
{}[x^{[2p]},y]=(\ad_x)^{2p}(y)\quad\text{~for any~}x\in\fg_\od,~y\in\fg.
\end{equation*}
The pair of maps $[p]$ and $[2p]$ is called a $p$-\textit{structure}
(or, sometimes, a $p|2p$-\textit{structure}) on $\fg$, and $\fg$ is
said to be \textit{restricted}. It suffices to determine the $p|2p$-structure on any basis of $\fg$; on simple Lie superalgebras there are not more than one $p|2p$-structure.

\sssbegin{Remark}\label{ResQuo} If (\ref{restr3}) is not satisfied,
the $p$-structure on $\fg_\ev$ does not have to generate a
$p|2p$-structure on $\fg$: even if the actions of $(\ad_x)^p$ and
$\ad_{x^{[p]}}$ coincide on $\fg_\ev$, they do not have to coincide
on the whole of $\fg$. For example, consider $p=2$
and 
$\fg=\foo^{(1)}_{I\Pi}(1|2)$ (for the definition of ortho-orthogonal Lie superalgebras, see \cite{LeP,
BGL2}; for clarity, we write $\fwk^{(1)}(3;a)$ or $\fg^{(1)}(A)$
instead of $(\fwk(3;a))^{(1)}$ or $(\fg(A))^{(1)}$ and the like) with basis $\{X_-^2, X_-, H, X_+, X_+^2\}$, where $X_-^2$,
$H$, and $X_+^2$ are even elements while $X_-$ and $X_+$ are odd
ones, with the relations (other bracket being equal to 0)
\[
{}[H,X_\pm]=X_\pm;\quad [X_+,X_-]=[X_+^2,X_-^2]=H.
\]
We can define a $2$-structure on $\fg_\ev\simeq\fsl(2)$, which is
nilpotent for $p=2$, by setting
\[
(X_-^2)^{[2]}=H;\quad H^{[2]}=H;\quad (X_+^2)^{[2]}=0,
\]
and extending it to the whole $\fg_\ev$ by properties
(\ref{restricted-1}) and (\ref{restricted-2}). This $2$-structure on
$\fg_\ev$ can not be extended to a $2|4$-structure on $\fg$, since,
for example,
\[
{}[X_-^2,[X_-^2,X_-]]=0\neq [(X_-^2)^{[2]}, X_-]=X_-.
\]
\end{Remark}

\sssec{Restricted modules}\label{resModP,2P} A $\fg$-module $M$
corresponding to a representation $\rho$ of the restricted Lie
superalgebra $\fg$ is said to be \textit{restricted} or having a
\textit{$p|2p$-structure} if
\begin{equation*} \label{restrP2P}
\begin{array}{ll}
\rho(x^{[p]})=(\rho(x))^{p}&\text{for any~}x \in\fg_\ev,\\
\rho(x^{[2p]})=(\rho(x))^{2p}&\text{for any~}x \in\fg_\od.
\end{array}\end{equation*}

\ssec{On $2|2$-structures on Lie superalgebras} Let $p=2$, a Lie
superalgebra $\fg$ have a $2|4$-structure, and $\textbf{F}(\fg)$
be the Lie algebra one gets from $\fg$ by forgetting the squaring
and considering only brackets by setting $[x,x]:=2x^2=0$ for $x$ odd.
Then $\textbf{F}(\fg)$ has a $2$-structure given by
\begin{equation}\label{22str}
\begin{array}{l}
\text{the ``2" part of $2|4$-structure on $\fg_\ev$};\\
\text{the squaring on $\fg_\od$, i.e., $x^{[2]}:=x^{2}$};\\
\text{the rule $(x+y)^{[2]}:=\begin{cases}x^{[2]}+y^{[2]}+[x,y]&\text{if
$x, y\in \fg_\ev$},\\
x^{2}+y^{[2]}+[x,y]&\text{if $x\in \fg_\od, y\in \fg_\ev$},\\
x^{2}+y^{2}+[x,y]&\text{if $x, y\in \fg_\od$}.\\
\end{cases}$}\\
\end{array}
\end{equation}
(Actually, the first and the third cases in \eqref{22str} are redundant. If $x$ and $y$ are both in $\fg_\ev$ or both in $\fg_\od$, then $x+y$ is homogenous, and $(x+y)^{[2]}$ in $\textbf{F}(\fg)$ is already given by $(x+y)^{[2]}$ or $(x+y)^{2}$, correspondingly.)
So one can say that if $p=2$, then the restricted Lie superalgebra
(i.e., the one with a $2|4$-structure) also has a
\textit{$2|2$-structure} which is defined even on inhomogeneous
elements (unlike $p|2p$-structures). In future, for Lie
superalgebras with $2|2$-structure, we write $x^{[2]}$ instead of
$x^{2}$ for any odd or inhomogeneous $x\in\fg$. The analog of sufficient condition 2) of Remarks~\ref{Uniq} holds.

\sssec{Restricted modules}\label{resMod22} A $\fg$-module $M$
corresponding to a representation $\rho$ of the Lie superalgebra
$\fg$ with $2|2$-structure is said to be \textit{restricted} or
having a \textit{$2|2$-structure} if
\begin{equation*} \label{restr22}
\begin{array}{ll}
\rho(x^{[2]})=(\rho(x))^{2}&\text{for any~}x \in\fg.
\end{array}\end{equation*}

\ssec{Restrictedness of Lie (super)algebras with Cartan matrix, and
of their relatives (the derived algebras, central extensions, and
quotients thereof modulo center)}{}~{}

\sssec{Lie (super)algebras with Cartan matrix} Speaking about Lie
(super)algebras $\fg=\fg(A)$ with an $n\times n$ Cartan matrix
$A=(A_{ij})$, recall (see \cite{LeD, BGL2}) that any nonzero element
$\alpha\in\Ree^n$ is called \textit{a root} if the homogeneous
subspace $\fg_\alpha$ of $\fg$ with grade (weight) $\alpha$ is nonzero. Let
$R$ be the set of all roots of $\fg$ and $\fh$ the maximal torus.

\parbegin{Proposition}\label{g(A)-p|2p} $1)$ If $p>2$ (or $p=2$ and
$A_{ii}\neq \od$ or $1$ for any $i$) and $\fg(A)$ is a Lie
(super)algebra, then $\fg(A)$ has a $p$-structure (resp.
$p|2p$-structure) such that
\begin{equation}\label{p-str}
\begin{array}{l}
(x_\alpha)^{[p]}=0\text{~for any even $\alpha\in R$ and $x_\alpha\in\fg_\alpha$},\\
(x_\alpha)^{[2p]}=0\text{~for any odd $\alpha\in R$ and $x_\alpha\in\fg_\alpha$},\\
\fh^{[p]}\subset \fh.
\end{array}
\end{equation}

$2)$ If $A_{ij}\in\Zee/p$ for all $i,j$, then the derived Lie (super)algebra
$\fg^{(1)}(A)$ inherits the $p$-structure (resp. $p|2p$-structure)
of $\fg(A)$, and we can make the $3$rd line of eq.~\eqref{p-str}
precise:
\begin{equation}\label{p-str11}
h_i^{[p]}=h_i\text{~~for any basis element $h_i\in\fh$}.
\end{equation}

$3)$ The quotient modulo center of a Lie (super)algebra $\fg$ with a
$p$-structure (resp. $p|2p$-structure) always inherits the
$p$-structure (resp. $p|2p$-structure) of $\fg$.
\end{Proposition}

\begin{proof}\label{rrr} 1) The statement of line 1 of eq.~\eqref{p-str} is
contained in the proof of \cite[v.1, Th.7.2.2]{S}, proof of the statement of line 2 is similar.

In our definition of roots, see \cite{BGL2}, instead of nonexisting
pairing $(\gamma, \, h)$, where $\gamma\in\Ree^n$ is a root and
$h\in\fh$, we introduce a function
\[
\evv: R\times\fh\to \Kee,\quad \evv(\gamma,\ h)=\text{the eigenvalue of
$\ad_h$ on $\fg_\gamma$.}
\]
Clearly, $\evv$ is linear in the second argument in the usual sense,
and it is $\Zee$-linear, i.e., additive, in the first argument in
the sense that
\[
\evv\left(\sum\limits_{1\leq i\leq k}
c_i\gamma_i,\ h\right)=\sum\limits_{1\leq i\leq k} c_i
\evv(\gamma_i,\ h) \text{~~for any $c_i\in \Zee$ and any roots
$\gamma_i$ such that $\sum c_i\gamma_i\in R$},
\]
here the $c_i$ are considered as elements of $\Zee$ or $\Ree$ in the
left-hand side and as elements of $\Kee$ in the right-hand side. The
other way round, if a function $f:R\to\Kee$ is such that
\begin{equation}\label{lin-fun}
f\left(\sum\limits_{1\leq i\leq k}
c_i\gamma_i\right)=\sum\limits_{1\leq i\leq k} c_i f(\gamma_i)\
\text{ for any $c_i\in \Zee$ and any roots $\gamma_i$ such that
$\sum c_i\gamma_i\in R$},
\end{equation}
then
\begin{equation}\label{why}
\text{there is an element $h\in\fh$ such that
$\evv(\gamma,h)=f(\gamma)$ for any root $\gamma$.}
\end{equation}

To prove the claim \eqref{why}, it suffices to choose $h$ so that the
condition $\evv(\alpha_j,h)=f(\gamma_j)$ is satisfied for any simple
positive root $\alpha_j$. Such an $h$ exists because the roots
$\alpha_j$ are linearly independent. Because any root can be
represented as an integral linear combination of simple (positive)
roots and due to the condition \eqref{lin-fun}, this condition is satisfied for
any root. This demonstrates validity of the statement of line 3 in
\eqref{p-str}.

Note that the derived Lie superalgebra $\fg^{(1)}(A)$ may fail to
have property \eqref{why}.

2) Note that for any integers $c_1,\dots,c_k$ and any roots
$\gamma_1,\dots,\gamma_k$ such that
$\gamma=\mathop{\sum}\limits_{1\leq i\leq k} c_i\gamma_i$ is a root,
the eigenvalue of $(\ad_h)^p$ on $\fg_\gamma$ is equal to
\[
\evv\left(\sum\limits_{1\leq i\leq k}
c_i\gamma_i,\ h\right)^p=\left(\sum\limits_{1\leq i\leq k} c_i
\evv(\gamma_i,\ h)\right)^p=\sum\limits_{1\leq i\leq k} (c_i
\evv(\gamma_i,\ h))^p=\sum\limits_{1\leq i\leq k} c_i
(\evv(\gamma_i,\ h))^p,
\]
since $c^p=c$ for any $c\in\Zee/p$, and $\binom{p}{m}\equiv 0 \mod
p$ for $m\not\equiv 0 \mod p$. In other words, the function $\evv$
on $\fh$ whose value at $h\in\fh$ is the eigenvalue of $(\ad_h)^p$
on $\fg_\gamma$, satisfies the condition \eqref{lin-fun}. Thus,
there really exists an element $h^{[p]}\in\fh$ such that
$\ad_{h^{[p]}}=(\ad_h)^p$.

If $A_{ij}\in\Zee/p$ for all $i,j$, we consider $A_{ij}$ as integers from the set $\{0, \dots, p-1\}$, and apply
Fermat's little theorem.

If $A_{ij}\not\in\Zee/p$ for some $i,j$, then $\fg^{(1)}(A)$ may have no
$p$-structure (resp. $p|2p$-struc\-ture) even if $\fg(A)$ has one.

3) Factorization modulo center fixes one $p$-structure (resp.
$p|2p$-structure) inherited from the nonfactorized algebra which may
have several such structures.
\end{proof}

\sssbegin{Examples}\label{EX} 1) Observe that the center $\fc$ of the Lie algebra
$\fwk(3; a)$ with Cartan matrix
$\footnotesize
A=\begin{pmatrix}
 \ev & a & 0 \\
a & \ev & 1 \\
 0 & 1 & \ev\end{pmatrix},
$
where $a\neq 0,1$, is spanned by $a h_1 +h_3$ in the standard numeration of Chevalley generators $x^{\pm}_i$ and
related elements $h_i:=[x^+_i, \ x^{-}_i]$. The
2-structure on $\fwk(3;a)$ is given by the conditions
$(x_\alpha)^{[2]}=0$ for any root vector $x_\alpha$, and the
conditions \eqref{2strwk3n}.

For the matrix $B=(1,0,0)$ supplementing the Cartan matrix $A$ as explained in \cite{BGL2} (where
in a copy of eq.~\eqref{2strwk3n} and similar ones there are typos:
$\ad_{h_i}^{[2]}$ and $\ad_{d}^{[2]}$ should be $h_i^{[2]}$ and
$d^{[2]}$, respectively) and the grading operator $d$ also defined in \cite{BGL2}, set:
\begin{equation}\label{2strwk3n}
\begin{array}{ll}
h_1^{[2]} =
(1+a t) h_1 + t h_3&\equiv h_1\pmod{\fc},\\
h_2^{[2]} = a t h_1 + a h_2 + t h_3 + (1 +a) d
&\equiv a h_2 + (1 +a) d\pmod{\fc},\\
h_3^{[2]} = (a t +a^2) h_1 + t h_3&\equiv a^2 h_1\pmod{\fc},\\
d^{[2]} = a t h_1 + t h_3 +d&\equiv d\pmod{\fc}.\\
\end{array}\end{equation}

Observe that the simple Lie algebra $\fwk^{(1)}(3; a)/\fc$ has no
2-structure.

2) The 2-structure on $\fwk(4;a)$ with Cartan matrix $\footnotesize\begin{pmatrix}
  \ev & a & 0& 0 \\
  a & \ev & 1 & 0 \\
  0 & 1 & \overline{0} & 1 \\
  0 & 0 & 1 & \ev
 \end{pmatrix},$  where $a\neq 0,1$, is given by
the conditions $(x_\alpha)^{[2]}=0$ for any root vectors $x_\alpha$,
and
\begin{equation}
\label{2strwk4}
\begin{array}{ll}
h_1^{[2]} =a h_1 +(1+a) h_4,&h_3^{[2]} = h_3,\\
h_2^{[2]} = a h_2,&h_4^{[2]} = h_4.\\
\end{array}
\end{equation}

\end{Examples}

\sssbegin{Proposition}\label{2vert4} The superizations of simple
Lie algebras of the form $\fg(A)$, and non-simple ones of the form $\fg(A)/\fc$, see \emph{\cite{BGL2}}, have $2|4$-structures.
\end{Proposition}

\begin{proof} The $2|4$-structure is given by means of the first two
lines of eq.~\eqref{p-str} and --- instead of the third line of eq.~
\eqref{p-str}
--- expressions \eqref{2strwk3n},
\eqref{2strwk4}, or \eqref{p-str11} if $A_{ij}\in\Zee/p$ for all
$i,j$.\end{proof}

\ssec{$(2,4)$- and $(2, -)$-structures on Lie algebras}\label{ss2,4}
Let $\fg=\fg_+\oplus \fg_-$ be a $\Zee/2$-grading of a Lie algebra
(not superalgebra) $\fg$. We say that $\fg$ has a \textit{$(2,
-)$-structure}, if there is a map \[
{}[2]:\fg_+\to\fg_+, \ x\mapsto x^{[2]}
\]
such that
(for simplicity, we consider the case of centerless $\fg$)
\[
{}[x^{[2]},y]=[x,[x,y]]\quad\text{~for any~}x\in\fg_+,\ \ y\in\fg,
\]
but there is no $2$-structure on $\fg$. 

A
 \textit{$(2,4)$-structure} (do not confuse with $2|4$-structure!) on a $\Zee/2$-graded Lie algebra $\fg$ is a pair: a $2$-structure on $\fg_+$, and  a map
\be\label{2,4}
{}[4]:\fg_-\to\fg_+,\ x\mapsto x^{[4]}\text{~~ such that $[x^{[4]},y]=[x,[x,[x,[x,y]]]]$~for any~}x\in\fg_-,\
y\in\fg.
\ee

The analog of sufficient condition 2) of Remarks~\ref{Uniq} holds. 

In the examples we know, the \textit{$(2,4)$-structure} is related with the following fact:
\begin{equation}\label{reason}
\begin{minipage}[c]{11 cm}
Whereas the maximal dimension of the irreducible $\fsl(2)$-module with
highest weight is equal to $p$ if $p>2$, the maximal
dimension of irreducible module with highest weight over the simple
3-dimensional Lie algebra for $p=2$ is equal to 4, see \cite{D}.
\end{minipage}
\end{equation}

\sssec{Examples}\label{Ex3} 1) If indecomposable symmetrizable $n\times
n$-matrix $A$ is such that
\[
A_{nn}=\od;\quad A_{ii}=\ev\text{~for~} i<n,
\]
then $\fg(A)$ has no $2$-structure but has a $(2,4)$-structure with
the $\Zee/2$-grading given on Chevalley generators by setting:
\[
\deg(\fh)=\ev;\quad \deg(x_n^\pm)=\od;\quad
\deg(x_i^\pm)=\ev\text{~~for~}i<n.
\]
In particular, the Lie algebra $\fg=\fo^{(1)}(2n+1)$ with Cartan
matrix
$\footnotesize
\begin{pmatrix}\ddots&\ddots&\ddots&\vdots\\\ddots
&\ev&1&0\\\ddots&1&\ev&1\\\cdots&0&1&\od\end{pmatrix}
$
can be viewed as the algebra of matrices of the form (recall that
$ZD(n)$ is the space of symmetric matrices with zeros on the main
diagonal)
\[\footnotesize
\begin{pmatrix}A&X&B\\Y^T&0&X^T\\C&Y&A^T\end{pmatrix},\text{~where~}\begin{array}{c}
A\in\fgl(n);~~B,C\in ZD(n);\\X,Y\text{~are
column-vectors.}\end{array}
\]
Then $\fg_+$ consists of matrices with $X=Y=0$ and the map ${}[2]$
is the squaring of matrices, while $\fg_-$ consists of matrices with
$A=B=C=0$, and the map ${}[4]$ is rising matrices to the fourth
power.

2) Let $p=(p_1, \dots, p_n)$, $q=(q_1, \dots, q_n)$, and one more indeterminate $x$ be
indeterminates of the generating functions of the Lie algebra
$\fh_I(2n+1;\underline{\One})$.
On it, there is a (2,4)-structure, see eq.~
\eqref{2,4}, with respect to the grading
\[
\text{$\deg p_i=\deg q_i=\ev$ for all $i$, $\deg x= \od$}.
\]

3) Example of a $\Zee/2$-graded Lie algebra $\fg$ which has a $(2, -)$-structure, but not a $(2, 4)$-structure: let $\fg$  be spanned by elements $a_n$, where $n\in\Zee$, and $b$. Let the commutation relations be 
\[
[a_m,a_n]=0, \ \ [b,a_n] = a_{n-1}\text{~~ for any $m,n\in\Zee$}.
\] 
Then $\fg$ has a $\Zee/2$-grading such that $\deg a_n = n\mod 2$ and $\deg b = \od$. This $\fg$ has a $(2, -)$-structure which maps every even element to $0$, but no element can be $b^{[4]}$.

\sssec{$(2,4)|4$-structure on Lie superalgebras} The Lie superalgebra $\fo\fo_{I\Pi}^{(1)}(2k_\ev+1|2k_\od)$
has, similarly to the
above,  a \textit{$(2,4)|4$-structure} consisting of the rising to the fourth
power
defined on the odd part while on the
 even part a $(2,4)$-structure is defined, for any
$y\in\fo\fo_{I\Pi}^{(1)}(2k_\ev+1|2k_\od)$,  satisfying the following conditions:
\begin{equation*}\label{2,4|2str}
\renewcommand{\arraystretch}{1.4}
\begin{array}{l}
\ad_{x^{[2]}}(y)=(\ad_x)^2(y)\text{~ for any $x\in
(\fo\fo_{I\Pi}^{(1)}(2k_\ev+1|2k_\od)_+)_\ev$}\\
\ad_{x^{[4]}}(y)=(\ad_x)^4(y)\text{~ for any $x\in
\fo\fo_{I\Pi}^{(1)}(2k_\ev+1|2k_\od)_-$ or $
(\fo\fo_{I\Pi}^{(1)}(2k_\ev+1|2k_\od)_+)_\od$.}
\end{array}
\end{equation*}

The analog of sufficient condition 2) of Remarks~\ref{Uniq} holds.

\sssec{$(2,4)|4$-restricted modules}\label{resMod2,4} A $\fg$-module
$M$ corresponding to a representation $\rho$ of the Lie (super)algebra
$\fg$ with $(2,4)$-structure (resp.
$(2,4)|4$-structure) is said to be \textit{$(2,4)$-restricted}
(resp. \textit{$(2,4)|4$-restricted}) or having a
\textit{$(2,4)$-structure} (resp. \textit{$(2,4)|4$-structure}) if
\begin{equation*} \label{restM2,4}
\begin{array}{ll}
\rho(x^{[2]})=(\rho(x))^{2}&\text{for any~}x \in(\fg_+)_\ev,\\
\rho(x^{[4]})=(\rho(x))^{4}&\text{for any $x \in\fg_-$ or $(\fg_+)_\od$}.
\end{array}\end{equation*}

\ssec{Restricted vectorial Lie
(super)algebras}\label{queervectorial} In Proposition \label{vect}
below, by ``vectorial Lie (super)al\-gebra" we only mean the one whose
the Weisfeiler filtration and grading corresponding to the maximal
subalgebra of least codimension (i.e., the degree of each
indeterminate is nonzero). For the list of known vectorial Lie
(super)algebras having simple derived, see \cite{GL3, Sk, BGL1,
BGLLS1, BGLLS2}.

Let $\fg(\sdim;\un)$ or $\fg(a;\un|b)$ denote the vectorial Lie
(super)algebra with a ``family name" $\fg$ realized on $\sdim=(a|b)$
indeterminates, of which $a$ are even and $b$ are odd, and with
shearing vector $\un$; if $b=0$ the notation usually shrinks to
$\fg(\dim;\un)$.

\sssec{Fact}\label{fact} Let vectorial Lie (super)algebra $\fg(\sdim;\underline{\One})$ be the
\textit{prolong}, i.e., the result of generalized (perhaps, partial)
Cartan prolongation, see \cite{Shch}, of the pair $(\fg_-, \fg_0)$, where $\fg_0$ is a
restricted Lie (super)algebra and $\fg_- =\mathop{\oplus}\limits_{-d\leq i<0}$ is a restricted $\fg_0$-module. Then $\fg(\sdim;\underline{\One})$ is
restricted. For the known simple derived of vectorial Lie
(super)algebras the restricted structure is given by the following
expressions, where $\fh$ is a maximal torus of $\fg_0$:
\begin{equation}\label{pvert2p}
\begin{minipage}[l]{12cm}
$\fh^{[p]}\subset \fh$, and if the structure constants lie in
$\Zee/p$, then\\ $h_i^{[p]}=h_i$ for
the basis elements $h_i$ of $\fh$;\\
$w^{[p]}=0$ (resp. $w^{[2p]}=0$) for the other even (resp. odd)
weight ele\-ments $w$ of the basis of $\fg$ with weights relative a maximal torus of $\fder \fg$.\end{minipage}
\end{equation}

\sssbegin{Proposition}\label{vect1} \emph{1)} For $p>0$,  the vectorial Lie superalgebra
$\fg(\sdim;\un)$ is not restricted if $\un\neq
\underline{\One}$.

\emph{2)} If $\fg:=\fg(\sdim;\One)$ is restricted, and the $i$-th derived $\fg^{(i)}$
of $\fg$ contains the maximal torus of $\fg$, then $\fg^{(i)}$ is
restricted and formulas \eqref{pvert2p} are applicable.
\end{Proposition}

\begin{proof} Items 1) and 2): for Lie algebra case, see \cite[v.1, Th. 7.2.2]{S};
the super case
is analogous. 
\end{proof}

\sssec{Remarks} 1) The restrictedness was established, so far, in
the following cases of simple Lie algebras: for $p>7$, by Block and
Wilson (\cite{BW}); for $p\geq 5$, see \cite[v.1, Th. 7.2.2]{S}; for the Lie
algebras with indecomposable Cartan matrix for any $p$, see \cite[Prop.
\ref{g(A)-p|2p}]{BGL2}.

2) Let $\bar x_i$ be the highest possible (divided) power of $x_i$,
and $\bar x=\prod \bar x_i$. Let
\[
\widetilde \fsvect(\dim;\un):=\{(1-\bar x)y \  \text{~~for any
$y\in\fsvect(\dim;\un)$}\}.
\]
The Lie algebra $\widetilde \fsvect(\dim;\underline{\One})$ is
restricted, but
expression \eqref{pvert2p} should be modified: the only difference
between bases of $\fsvect(\dim;\underline{\One})$ and $\widetilde \fsvect(\dim;\underline{\One})$ is that instead of $\partial_i$ we
have $X_i:=(1-\bar x)\partial_i$ and while
$\partial_i^{[p]}=0$ we have $X_i^{[p]}= - (\del_i^{p-1} \bar{x}) \partial_i$, see \cite{BKLLS}.


\section{The two methods of superization for $p=2$}\label{Sproof}

\ssec{Queerification for $p\neq 2$}\label{111} Let $\cA$ be an
associative algebra; let $\cA_L$ be the Lie algebra with the same space as $\cA$ and
the multiplication being defined by the commutator instead of the dot product.
The space of the Lie superalgebra $\fq(\cA)$, which we call the
\textit{queerification} of $\cA$, is $\cA_L\oplus\Pi(\cA)$, where
$\Pi$ is the change of parity functor,
so $\fq(\cA)_\ev=\cA_L$ and $\fq(\cA)_\od=\Pi(\cA)$, with the multiplication given
by the following expressions
\begin{equation}\label{commRel}
{}[x,y]=xy-yx;\quad [x,\Pi(y)]=\Pi(xy-yx);\quad [\Pi(x),\Pi(y)]=
xy+yx \ \text{~for any~~}x,y\in \cA.
\end{equation}
The term ``queer", now classical, is taken after the Lie
superalgebra $\fq(n):=\fq(\Mat(n))$, a \lq\lq queer\rq\rq\ analog (as explaned in \cite{Lsos}) of $\fgl(n)$, where $\Mat(n)$ is the
associative algebra of $n\times n$ matrices. We express the elements
of the Lie superalgebra $\fg=\fq(n)$ by means of a pair of matrices
\begin{equation}
\label{.1} (A,B)\longleftrightarrow
\begin{pmatrix}A&B\\B&A
\end{pmatrix}
\in \fq(n), \text{~~where $A,B \in \fgl(n)$}.
\end{equation}

We will similarly denote by pairs $(A,B)$, where $A,B \in\cA$ the elements of $\fq(\cA)$. The brackets between these basis elements are as follows:
\begin{equation}\label{.2}
\renewcommand{\arraystretch}{1.4}
\begin{array}{l}
{}[(A_1,0),(A_2,0)]=([A_1,A_2],0),\quad
[(A,0),(0,B)]=(0,[A,B]),\\
{}[(0,B_1), (0,B_2)]=(B_1B_2+B_2B_1,0).
\end{array}\end{equation}

\sssec{The only simple Lie superalgebras related with queerification
for $p\ne2$ are $\fp\fs\fq(n)$ for $n>2$} Let $\fs\fq(n):=\fq(n)^{(1)}$ denote the
subsuperalgebra of queertraceless matrices, where the \textit{queer
trace} is $\qtr:(A,B) \mapsto \tr B$. The Lie superalgebras $\fq(n)$ and
$\fs\fq(n)$ are specifically ``super" analogs of the general Lie
algebra $\fgl(n)$ and its special (traceless) subalgebra $\fsl(n)$; we define their projectivizations to be $\fp\fq(n):=\fq(n)/\Kee 1_{2n}$ and
$\fp\fs\fq(n):=\fsq(n)/\Kee 1_{2n}$.

\ssec{Queerification for $p=2$}\label{112} For $\fg=\fq(\cA)$, where $\cA$ is an associative algebra, the
multiplication is defined by the expressions \eqref{.2}, the bracket
of odd elements being polarization of squaring of the odd elements:
\begin{equation*}\label{.2for2}
(0,B)^2=(B^2, 0).
\end{equation*}

If $p=2$, then we
can queerify any restricted Lie algebra $\fg$, not only associative algebras, as follows. We set
$\fq(\fg)_\ev=\fg$ and $\fq(\fg)_\od=\Pi(\fg)$; define the
multiplication involving the odd elements as follows:
\begin{equation}\label{q(g)}
{}[x,\Pi(y)]=\Pi([x,y]);\quad (\Pi(x))^2=x^{[2]}\quad\text{~for
any~~}x,y\in\fg.
\end{equation}
Clearly, if $\fg$ is restricted and $\fii\subset\fq(\fg)$ is an ideal,
then $\fii_\ev$ and $\Pi(\fii_\od)$ are ideals in $\fg$. So, if $\fg$ is
restricted and simple, then $\fq(\fg)$ is a simple Lie superalgebra.
(Note that $\fg$ has to be simple as a Lie algebra, not just as a
\textit{restricted} Lie algebra, i.e., $\fg$ is not allowed to have
\textbf{any} ideals, not only restricted ones.) A generalization of
the queerification is the following procedure producing as many
simple Lie superalgebras as there are simple Lie algebras.

\sssec{\underline{Method 1}: \textbf{generalized queerification}}\label{Method1} Let the \textit{$1$-step
restricted closure} $\fg^{<1>}$ of the simple Lie algebra $\fg$ be
the minimal subalgebra of the (classically) restricted closure
$\overline{\fg}$ containing $\fg$ and all the elements $x^{[2]}$,
where $x\in\fg$. To any simple Lie algebra $\fg$ the \textit{generalized queerification} assigns the Lie superalgebra
\begin{equation*}\label{tildeQ}
\tilde \fq(\fg):=\fg^{<1>}\oplus \Pi(\fg)
\end{equation*}
with squaring given by $(\Pi(x))^2=x^{[2]}$ for any $x\in\fg$.
Obviously, for $\fg$ restricted, the generalized queerification coincides with the
queerification: $\tilde \fq(\fg)=\fq(\fg)$.

\sssbegin{Theorem}\label{MainTh} For  any simple Lie algebra $\fg$, the Lie superalgebra
$\tilde \fq(\fg)$ is simple.
\end{Theorem}

\begin{proof} Assume that $\fii = \fii_\ev\oplus
\Pi(\fii_\od)$, where $\fii_\ev\subset \fg^{<1>}$ and $\fii_\od\subset \fg$,
is an ideal of $\tilde \fq(\fg)$. Then $\fii_\od$ is an ideal in $\fg$
since $[\fg, \Pi(\fii_\od)]\subset \Pi(\fii_\od)$, so either $\fii_\od=0$, or
$\fii_\od=\fg$. If $\fii_\od=\fg$, then, by construction, $\Pi(\fii_\od)$
generates the whole $\fg^{<1>}$, and $\fii=\tilde \fq(\fg)$, so in what
follows we assume that $\fii_\od=0$.

Similarly, $\fii_\ev$ is an ideal of $\fg^{<1>}$. Since $\fg$ is an
ideal in $\overline{\fg}$, it follows that $\fii_\ev\cap\fg$ is an
ideal of $\fg$, so either $\fii_\ev\cap\fg=0$, or $\fg\subset \fii_\ev$.
If $\fg\subset \fii_\ev$, then $\fii_\od\ne 0$ and this is a contradiction
with our assumption.

If $\fii_\ev\cap\fg=0$, then $\fii_\ev$ commutes with $\fg$. It follows
from the minimality of the restricted closure that the centralizer
of $\fg$ in $\overline{\fg}$ is the center of $\fg$, so in this case
$\fii_\ev=0$, and $\fii=0$. \end{proof}

For a detailed description of several examples of generalized queerifications, see \S\ref{Sq}. Additionally, let us also mention
$\tilde \fq(\fwk^{(1)}(3;a)/\fc)$, and numerous new simple Lie
superalgebras that are generalized queerifications of simple vectorial
Lie algebras (of serial type, such as $\fvect$, $\fsvect$, and
$\fh_I$, $\fh_\Pi$, $\fk$ with their divergence-free subalgebras,
see \cite{LeP,ILL}, and their deforms, see \cite{Sk, Kos}; and
exceptional ones, see \cite{BGLLS,BGLLS1, BGLLS2}).

\ssec{\underline{Method 2}: superization of the $\Zee/2$-graded
simple Lie algebra $\fg$}\label{Method2} Observe that every simple Lie algebra has a $\Zee/2$-grading because every simple Lie algebra $\fg$ has a nonzero toral subalgebra.

For $p=2$, the space of \textbf{every} $\Zee/2$-graded
simple Lie algebra $\fg$ can be endowed with a Lie superalgebra structure (more precisely,
unless $\fg$ is restricted, its space has to be enlarged prior to
superization, as described below). For restricted Lie algebras $\fg$, these ``hidden supersymmetries" were first
discovered in \cite{BGL2}, but their
mechanism remained unclear until
now.

If $p=2$, let $\fg=\fg_+\oplus\fg_-$ be a simple Lie algebra with a
$\Zee/2$-grading $\gr$. Let $(\fg,\gr)$ be the minimal Lie subalgebra of
the restricted closure $\overline{\fg}$ containing $\fg$ and all the
elements $x^{[2]}$, where $x\in\fg_-$. Clearly, there is a single
way to extend the grading $\gr$ from $\fg$ to $(\fg,\gr)$. 

Let $\fs(\fg,\gr)$ be
the Lie superalgebra structure on the space of $(\fg,\gr)$ given by
\[
x^2:=x^{[2]}\text{~~for any $x\in\fg_-$.}
\]

\sssbegin{Theorem}\label{Statement} For  any simple Lie algebra $\fg$, the Lie superalgebra
$\fs(\fg,\gr)$ is simple. \end{Theorem}

\begin{proof} Let $\fii$ be an ideal of
$\fs(\fg,\gr)$. Let $\textbf{F}: \fs(\fg,\gr)\tto (\fg,\gr)$ be the
desuperization functor. Then $\textbf{F}(\fii)$ is an ideal in
$(\fg,\gr)$. Since $\fg$ is an ideal in $\overline{\fg}$, and
therefore in $(\fg,\gr)$, we see that $\textbf{F}(\fii)\cap\fg$ is an
ideal in $(\fg,\gr)$, and therefore in $\fg$. This means that either
$\fg\subset \textbf{F}(\fii)$, or $\textbf{F}(\fii)\cap\fg=0$. By
construction, $\textbf{F}^{-1}(\fg)$ generates $\fs(\fg,\gr)$, so if
$\fg\subset \textbf{F}(\fii)$, then $\fii=\fs(\fg,\gr)$. If
$\textbf{F}(\fii)\cap\fg=0$, then $\textbf{F}(\fii)$ commutes with $\fg$;
it follows from the minimality of the $2$-closure that the
centralizer of $\fg$ in $\overline{\fg}$ is the center of $\fg$, so
in this case $\fii=0$.
\end{proof}

For a huge number of new simple
Lie superalgebras, see \S\ref{Z2gra}, examples in
\cite{BGLLS1,BGLLS2} and references therein, and \cite{KrLe}.

\ssec{A relation between queerification and \lq\lq method 2"}
Let $\fg$ be an arbitrary Lie algebra over any field $\Kee$, and $C$ an associative and commutative algebra over $\Kee$. On the space $\fg_C:=\fg\otimes_{\Kee} C$, a Lie algebra structure is naturally defined (\lq\lq extension of the base field").

In particular, if $A=\Kee[a]/(a^2-1)$, then $\fg_A=\fg\oplus \fg a$, as spaces. Even if $\fg$ is simple, the Lie algebra $\fg_A$ is never simple, since the subspaces $\fg\otimes (1\pm a)$ are ideals in it. Assuming that $a$ is odd, we see that $\fg_A$ has a natural $\Zee/2$-grading.

Now, let $\Char\Kee=2$  and $\fg$ be simple. Then we can apply \lq\lq method 2" to $\fg_A$ considered with the natural grading. Then $\fs(\fg_A)$ is a simple Lie superalgebra.  This $\fs(\fg_A)$ is, as is easy to see, a queerification of $\fg$.

\section{Classification of simple Lie algebras yields the
classification of simple Lie superalgebras if $p=2$}\label{Snew}

Let $p=2$, and $\fg=\fg_\ev\oplus \fg_\od$ a Lie superalgebra,
$S=\Span\{x^2\mid x\in\fg_\od\}$.

\ssbegin{Lemma}\label{L1} The space $([\fg_\od,\fg_\od]+S)\oplus
\fg_\od$ is an ideal of $\fg$. \end{Lemma}

\begin{proof} The subspace $[\fg_\od,\fg_\od]$
is $\fg_\ev$-invariant due to Jacobi identity.
We see that
\[
{}[x^2,y]=[x,[x,y]]\in
[\fg_\od,\fg_\od]\text{~~for any $x^2\in S,\  y\in\fg_\ev$.}
\]
This means that $([\fg_\od,\fg_\od]+S)\oplus
\fg_\od$ is an ideal of $\fg$. \end{proof}

Let now $\fg$ be a simple Lie superalgebra. By Lemma \ref{L1}, we
see that $\fg_\ev=[\fg_\od,\fg_\od]+S$. Then
\begin{equation}\label{algFH}
\fh:=[\fg_\od,\fg_\od]\oplus \Pi(\fg_\od) \end{equation}
is an ideal of the Lie
algebra $F(\fg)$, whereas $\fg$ is obtained from $\fh$ by means of
``method 2". Notice that the $\Zee/2$-grading of $\fg$ induces $\Zee/2$-gradings on the Lie algebras $F(\fg)$ and $\fh$. In particular,  
\[
\text{$\fh=\fh_+\oplus \fh_-$, where $\fh_+=[\fg_\od,\fg_\od]$ and $\fh_-=\Pi(\fg_\od)$.}
\]

A question arises: must the Lie algebra $\fh$ be
simple?

\ssbegin{Lemma}\label{L2} The $\Zee/2$-graded Lie algebra $\fh$ defined by eq.~\eqref{algFH} for a simple $\fg$ has no nontrivial $\Zee/2$-graded
ideals. \end{Lemma}

\begin{proof} Let $\fii=\fii_+\oplus \fii_-$ be a $\Zee/2$-graded ideal of $\fh$.
This means that $[\Pi(\fg_\od),
\fii]\subset \fii$, and hence $[S,\fii]\subset \fii$, implying that $\fii$
is an ideal of the whole Lie algebra $F(\fg)$. Then
\[
I=\fii+\Span\{x^2\mid x\in\fii_-\}\subset F(\fg),
\]
is an ideal that can be superized by means of ``method 2".

Indeed, if $x\in \fii_-$, $y\in F(\fg)$, then
$[x^2,y]=[x,[x,y]]\in\fii$. If now $\fs(I)$ is the superization of
$I$, then $\fs(I)$ is an ideal of the Lie superalgebra $\fg$, and
since $\fg$ is simple, it follows that either $\fs(I)= 0$, or
$\fs(I)=\fg$. In the first case, $\fii=0$, while in the second case
$\fii\supset \Pi(\fg_\od)$, and hence, $\fii=\fh$. \end{proof}

What nongraded ideal $\fii$ might the Lie algebra $\fh$ have?

Denote the projectors of $\fh$ onto $\fh_+$ and $\fh_-$ by
$\pr_+$ and $\pr_-$, respectively. Then $(\fii\cap
\fh_+)\oplus (\fii\cap \fh_-)$ and $\pr_+(\fii)\oplus
\pr_-(\fii)$ are graded ideals of the Lie algebra $\fh$. By Lemma
\ref{L2} we see that
\[
(\fii\cap \fh_+)=0, \ \ (\fii\cap \fh_-)=0, \ \
\pr_+(\fii)=\fh_+, \ \ \pr_-(\fii)=\fh_-.
\]
This means that there exists a bijection $f:\fh_+\longrightarrow
\fh_-$ such that $\fii=\{x+f(x)\mid x\in \fh_+\}$. Now, since
$\fii$ is an ideal of $\fh$, it immediately follows that
\begin{equation}\label{eq1}
[y,f(x)]=f([y,x]), \quad
[f(y),f(x)]=f^{-1}([f(y),x])=[y,x] \text{~~for any
$x,y\in\fh_+$.}
\end{equation}
Eq. \eqref{eq1} means that $\fg=\tilde \fq(\fh_+)$.

\ssbegin{Theorem}\label{MainTh1} If $p=2$, then every simple
finite-dimensional Lie superalgebra can be obtained from a simple
Lie algebra by means of either (generalized) queerification, or ``method
$2$". \end{Theorem}

\begin{proof} As it is written after Lemma~\ref{L1}, we can obtain $\fg$ by ``method 2"
from the Lie algebra $\fh$, see eq.~\eqref{algFH}. Just after Lemma~\ref{L2} it is shown
that if $\fh$ is not simple, it may only have a nontrivial \textit{nongraded} ideal, and
then $\fg = \tilde \fq(\fh_+)$. And the Lie algebra $\fh_+$ must be simple,
because, if $\fh_+$ had a nontrivial ideal $\fii$, then $\fg$
would have had the nontrivial ideal $\fg_\ev\oplus \Pi(\fii)$.

So if the Lie algebra $\fh$ is simple, then $\fg$ can be obtained by
``method 2" from $\fh$. If $\fh$ is not simple, then $\fg$ can
be obtained by queerification of the simple Lie algebra
$\fh_+$.\end{proof}

\section{Examples of queerifications for $p=2$}\label{Sq}

\ssec{Queerification of superalgebras}\label{QSup} If $\cB$ is an
associative \textbf{super}algebra, we can construct $\fq(\cB)$
using eq.~\eqref{commRel}, where $\cA=\cB_\ev\oplus \Pi(\cB_\od)$.  If
$p=2$, then we can also queerify any restricted Lie superalgebra $\fG$ using
eq.~\eqref{q(g)} for $\fg=\fG_\ev\oplus \Pi(\fG_\od)$.

For the definition of the $2|4$-structure on Lie superalgebras --- an analog of the 2-structure on Lie algebras,
see Subsection~\ref{SSp2pStr}

\sssbegin{Lemma}\label{queerexcept} Let  $p=2$.

$1)$ Let $\fg$ be a simple Lie
algebra with a $2$-structure. Every Lie superalgebra of the form
$\fq(\fg)$ has a natural $2|4$-structure given by the $2$-structure
on $\fg$.

$2)$ If $\cB$ is an associative superalgebra, then
 $\fq(\cB)\simeq \fq(\textbf{F}(\cB))$. If $\fg$ is a Lie superalgebra with a
$2|4$-structure, then $\fq(\fg)\simeq \fq(\textbf{F}(\fg))$.
\end{Lemma}

\begin{proof} Because $\fg$ is restricted, it follows that for every $x, y\in \fg$,
we have
\[
[x^{[2]},\Pi(y)]=\Pi[x^{[2]},y]=\Pi(\ad_x^2 (y))=\Pi([x,[x,y]])=
[x,[x,\Pi(y)]]=\ad_x^2 (\Pi(y)).
\]
Since condition \eqref{restr3} is met, the $2|4$-structure can be
defined as in subsec. \ref{SSp2pStr}.

Heading 1) follows from eq.~\eqref{commRel}; to prove
heading 2) recall eq.~\eqref{q(g)} as well.
\end{proof}

\ssec{Queerifications of $\fsl(n)$ and $\fpsl(n)$}\label{3.2.1}
If $p=2$, the supertrace on the Lie superalgebra $\fq(n)$ does not vanish identically as it does for $p\neq 2$, so there are two traces on $\fq(n)$ which induce two NISes (nondegenerate invariant symmetric bilinear forms) on the simple subquotient of $\fq(n)$ which is a \textit{double extension} of this simple subquotient, see \cite{BKLS}. This simple subquotient of $\fq(n)$ is a queerification of $\fsl(n)$ if $n$ is odd and $\fpsl(n)$ if $n$ is even.

\ssec{Queerification of orthogonal Lie algebra
$\fo_B(n)$}\label{q(o_B)} The Lie algebra $\fo_B(n)$ preserving a
nondegenerate symmetric bilinear form $B$ on the $n$-dimensional space has a
natural $2$-structure: considered as Lie algebra of $n\times n$ matrices (or linear
operators) $X$ such that $XB+BX^T=0$, we see that
\[
X^2B+B(X^2)^T=X(XB+BX^T)+(XB+BX^T)X^T=0.
\]
Thus, if $X\in\fo_B(n)$, then $X^2\in\fo_B(n)$. So we can consider
queerification $\fq(\fo_B(n))$ of this algebra. Recall (see \cite{LeP}) that there is one
equivalence class of bilinear forms $B$ if $n$ is odd, and two classes if $n$ is even. For the normal
shape of the Gram matrix of the form of one of the two classes one can take $B=I$, the unit matrix, for the normal shape of the Gram matrix  of the form of the other class we take either $B=\Pi$, the matrix with units on the side diagonal, the other matrix
elements being 0, or an equivalent to it $B=\footnotesize\begin{pmatrix} 0_n&1_n\\1_n&0_n\end{pmatrix}$.
We will find
simple (modulo center) derived of these queerifications for $n$
large enough; more precisely, $n\geq 3$ for $\fq(\fo_I(n))$, and
$n\geq 6$ for $\fq(\fo_\Pi(n))$.

First, consider the ``infinite derived"
\[
\fq(\fo_B(n))^{(\infty)}:=\bigcap\limits_{i\geq 1}
\fq(\fo_B(n))^{(i)}.
\]
Clearly, $\fq(\fo_B(n))^{(\infty)}_\ev$ and
$\Pi(\fq(\fo_B(n))^{(\infty)}_\od)$ are subalgebras of $\fo_B(n)$.

Let us give a general description of nontrivial ideals of
$\fq(\fo_B(n))^{(\infty)}$. Clearly, for any ideal $\fii\subset
\fq(\fo_B(n))^{(\infty)}$, we have
\begin{equation*}\label{wehave}
\begin{array}{ll}
1)&\fii_\ev\text{~~is an ideal of
$\fq(\fo_B(n))^{(\infty)}_\ev$;}\\
2)&\Pi(\fii_\od)\text{~~is an ideal of
$\Pi(\fq(\fo_B(n))^{(\infty)}_\od)$.}
\end{array}
\end{equation*}

Let us specify our description. We'll need a new notation $\widehat\fq$ for the quotient of $\tilde\fq$ or its derived:

\sssec{$\widehat\fq(\fo_I(n))$}\label{fq(fo_I(n))} Recall that $\fo_I(n)$
consists  of all symmetric matrices. Let $ZD(n)$ denote the space of symmetric $n\times n$
matrices with zeros on the main diagonal. We have
\[
\fo_I(n)^{(\infty)}=\fo_I(n)^{(1)}=ZD(n).
\]
It follows that
\[
\begin{array}{ll}
\fq(\fo_I(n))^{(1)}=\fo_I(n)\oplus \Pi(ZD(n));&\\
\fq(\fo_I(n))^{(i)}=(\fo_I(n)\cap\fsl(n))\oplus
\Pi(ZD(n))\text{~~for $i>1$}.
\end{array}
\]

Clearly, the Lie algebra  defined on the space $ZD(n)$ is
simple, so any nontrivial ideal of $\fq(\fo_I(n))^{(\infty)}$ has
zero odd part. The restricted Lie algebra $\fo_I(n)\cap\fsl(n)$ has
a nontrivial ideal if and only if $n$ is even, and this ideal is
$\Kee 1_n$. Thus,
\[
\widehat \fq(\fo_I(n)):=\begin{cases}(\fo_I(n)\cap\fsl(n))\oplus\Pi(ZD(n))&\text{if
$n$ is
odd},\\
(\fo_I(n)\cap\fsl(n))/\Kee1_n\oplus\Pi(ZD(n))&\text{if $n$ is even.}
\end{cases}
\]

\sssec{$\widehat\fq(\fo_\Pi(2k))$}\label{fq(fo_Pi(n))} The Lie algebra
$\fo_\Pi(2k)$ consists of matrices of the following form:
\[
\begin{pmatrix}
A&C\\D&A^T\end{pmatrix},\quad\text{where~}
\begin{array}{l}A\in\fgl(k),\\
\text{$C=C^T$ and $D=D^T$.}\end{array}
\]
So the elements of $\fq(\fo_\Pi(2k))$ are of the form
\[
\begin{pmatrix}
A&C\\D&A^T\end{pmatrix}\oplus\Pi\left(\begin{pmatrix}
A'&C'\\D'&A'^T\end{pmatrix}\right),\quad\text{where~}
\begin{array}{l}A,A'\in\fgl(k),\\\text{$C,C',D,D'$ are
symmetric matrices.}\end{array}
\]
Computations show that elements of $\fq(\fo_\Pi(2k))^{(i)}$ satisfy the
following conditions for different $i$:
\[
\begin{array}{ll}
i=1: &A,A'\in\fgl(k),\quad C,C',D,D'\in ZD(k);\\
i=2: &A\in\fgl(k),\quad A'\in\fsl(k),\quad C,C',D,D'\in ZD(k);\\
i\geq 3: &A,A'\in\fsl(k),\quad C,C',D,D'\in ZD(k). \end{array}
\]
Thus,
\[
\fq(\fo_\Pi(2k))^{(\infty)}_\ev\simeq\Pi(\fq(\fo_\Pi(2k))^{(\infty)}_\od)=\fo_\Pi(2k)^{(2)}=
\fo_\Pi(2k)^{(\infty)}.
\]
Clearly, $\fo_\Pi(2k)^{(2)}$ has a nontrivial ideal if and only if
$k$ is even, and this ideal is $\Kee 1_{2k}$, so
\[
\widehat \fq(\fo_\Pi(2k)):=\begin{cases}\fo_\Pi(2k)^{(2)}\oplus\Pi(\fo_\Pi(2k)^{(2)})&\text{if
$k$ is
odd},\\
(\fo_\Pi(2k)^{(2)}\oplus\Pi(\fo_\Pi(2k)^{(2)}))/(\Kee1_{2k}\oplus
\Kee\Pi(1_{2k}))&\text{if
 $k$ is even.}
\end{cases}
\]

\section{A relation between Cartan prolongation and queerification}\label{SrelCartAndQ}

In what follows, we illustrate the queerification phenomenon in
terms of Cartan prolongations for the simplest example. This illustration shows why the characteristic
$p=2$ is exceptional. For a detailed description of the generalized
Cartan prolongation and recipes how to realize elements of Lie
(super)algebras by vector fields, see \cite{Shch}. In this section, the results do not depend on
the shearing vector $\un$ (for related definitions, see \cite{BGLLS2}), so we mostly do not indicate it.

\ssec{Notation} Recall that we denote the
elements of $\fg:=\fq(n+1)$ by pairs of matrices $(A,B)$, where $A,B \in
\fgl(n+1)$, see eq.~\eqref{.1}. Let us shorthand the element
$(E_{ij},0)\in\fg_\ev$, see eq.~\eqref{.1}, by means of the matrix unit
$E_{ij}$, and $(0,E_{ij}) \in\fg_\od$ by means of the same matrix
unit but denoted $X_{ij}$ and considered to be odd. Let 0 through
$n$ be labels of rows and columns of matrices $A,B \in \fgl(n+1)$.
Consider the following $\Zee$-grading of $\fg=\fg_{-1}\oplus
\fg_0\oplus \fg_1$:
\begin{equation}\label{.11}
\renewcommand{\arraystretch}{1.4}
\begin{array}{l}
\fg_{-1}=\Span(E_{0i},\ X_{0i}\mid i=1,\dots, n),\quad
\fg_{1}=\Span(E_{i0},\ X_{i0}\mid i=1,\dots, n),\\
\fg_0=\Span(E_{00},\ X_{00},\ E_{ij},\ X_{ij}\mid i,j=1,\dots,
n)\cong \fq(n)\oplus \fq(1), \text{~~where}\\
\text{$\fq(n)=\Span(E_{ij},\ X_{ij}\mid i,j=1,\dots, n)$ and
$\fq(1)=\Span(E_{00},\ X_{00})$.}\end{array}
\end{equation}

\sssec{Homomorphism $\varphi:\fq(n+1)\tto \fvect(n|n)$} We define
$\varphi$ by setting
\begin{equation}\label{.12}
\renewcommand{\arraystretch}{1.4}
\begin{array}{l}
E_{0i}\mapsto \del_{z_i}, \quad X_{0i}\mapsto
\del_{\xi_i};\\

\begin{array}{lll}
E_{ij}\mapsto z_i\del_{z_j}+\xi_i\del_{\xi_j}, & X_{ij}\mapsto
z_i\del_{\xi_j}+\xi_i\del_{z_j}&\text{for
$i,j>0$,}\\
E_{00}\mapsto
-\mathop{\sum}\limits_i\left(z_i\del_{z_i}+\xi_i\del_{\xi_i}\right),
& X_{00}\mapsto\mathop{\sum}\limits_i\left(-
z_i\del_{\xi_i}+\xi_i\del_{z_i}\right);&\\
\end{array}\\

\begin{array}{l}
E_{i0}\longmapsto -z_i\mathop{\sum}\limits_j\left(z_j\del_{z_j}+
\xi_j\del_{\xi_j}\right)-\xi_i\mathop{\sum}\limits_j\left(
z_j\del_{\xi_j}-\xi_j\del_{z_j}\right),\\
X_{i0}\longmapsto z_i\mathop{\sum}\limits_j\left(-z_j\del_{\xi_j}+
\xi_j\del_{z_j}\right)-\xi_i\mathop{\sum}\limits_j\left(
z_j\del_{z_j}+\xi_j\del_{\xi_j}\right).
\end{array}\\
\end{array}
\end{equation}
The homomorphism $\varphi:\fq(n+1)\tto \fvect(n|n)$ has a kernel
yielding an embedding
\[
\varphi:\begin{array}{ll}
\fp\fq(n+1)\tto\fvect(n|n)&\text{if $p\ne 2$},\\
\fp\fs\fq(n+1)\tto\fvect(n|n)&\text{if $p= 2$.} \end{array}
\]
Let $J$ be  the odd operator
commuting with $\varphi(\fq(n))$, and $\varphi(\fg_{-1})$ the tautological
$\varphi(\fq(n))$-module.  We have
$\varphi(\fg_0)\cong \begin{cases}\fq(n)&\text{if $p= 2$,}\\
\fq(n)\oplus \Kee\cdot
J&\text{if $p\neq 2$.}\end{cases}$

\sssec{Prolongs of the nonpositive part of $\fp\fq(2)$ for $p\ne 2$}
In this particular case eqs.~\eqref{.12} take the form:
\begin{equation}\label{.101}\footnotesize
\renewcommand{\arraystretch}{1.4}
\begin{array}{llll}
E_{01}&\mapsto \pd z, &X_{01}&\mapsto \pd \xi,
\\
E_{00}&\mapsto -z\pd z-\xi\pd \xi, &
X_{00}&\mapsto -z\pd \xi+\xi\pd z, \\
E_{11}&\mapsto z\pd z+\xi\pd \xi,& X_{11}&\mapsto z\pd \xi+\xi\pd z, \\
E_{10}&\mapsto -z^2\pd z-2z\xi\pd \xi, &X_{10}&\mapsto -z^2\pd \xi.
\end{array}\end{equation}
\normalsize Clearly, $X_{11}-X_{00}\mapsto 2z\pd\xi$ and
$X_{11}+X_{00}\mapsto 2\xi\pd z$, so
$\varphi(\fg_0)\simeq\fsl(1|1)$, and hence
\[
(\fg_{-1},\varphi(\fg_0))_{*,\un}\simeq\fsvect(1;\un|1),
\]
the Lie
superalgebra of divergence-free vector fields, see \cite{BGLLS2}.

Denote the 0th component $\varphi(\fg_0)$ of the image of
$\fp\fs\fq(2)$ by $\fh$. Direct computations show:
\[
\begin{array}{l}
\fh=\Span\{z\pd z+\xi\pd\xi,
z\pd\xi\};\\
(\varphi(\fg_{-1}), \fh)_{*,\un} =\Span\{f(z)\pd
z+f'(z)\xi\pd\xi, \ \ g(z)\pd\xi\mid f,g\in \cO(1;\un)\}.
\end{array}
\]

Consider the \textit{generalized Cartan prolong} $\varphi(\fq(2))_{*,\un}$, i.e.,
the maximal Lie subsuperalgebra of $\fvect(1;\un|1)$ whose $(-1)$-st, $0$-th and
$1$-st components coincide with the respective components of
$\varphi(\fq(2))$. It is easy to see that this prolong coincides
with the semi-direct sum (here $\fii\ltimes\fa$ is a direct sum as spaces where $\fii$ is an ideal)
\[
\varphi(\fq(2))_{*,\un}= (\varphi(\fg_{-1}, \fh))_{*,\un} \ltimes \Kee\xi\pd z.
\]

An interpretation: the Cartan prolong $(\varphi(\fg_{-1}), \fh)_{*,\un}$
consists of the divergence-free vector fields preserving the
subspace of $\varphi(\fg_{-1})$ spanned by $\pd\xi$. Indeed, the
subalgebra $\fh\subset\fsl(1|1)$ reducibly acts on the component
$\varphi(\fg_{-1})$ preserving $\Kee\pd\xi$.

\sssec{Prolong of the nonpositive part of $\fp\fq(2)$ for $p=2$} For
$p=2$, all formulas of eq.~\eqref{.101} hold but the element
$X_{11}-X_{00}$ belongs to the kernel of the homomorphism, and hence
\begin{equation}\label{p=2}
\begin{array}{ll}
{}[E_{10}, E_{01}]=E_{11}-E_{00}\mapsto 0, & [E_{10},
X_{01}]=X_{11}-X_{00}\mapsto 0,\\
{} [X_{10}, E_{01}]=X_{11}-X_{00}\mapsto 0, & [X_{10},
X_{01}]=X_{11}+X_{00}\mapsto 0.
\end{array}
\end{equation} In other words, the images of $E_{10}$ and $X_{10}$ in
$\fsvect(1|1)$ must commute with every element of
$\varphi(\fg_{-1})$, hence these images should vanish. So the map
$\varphi$ has a big kernel and does not define embedding of either $\fq(2)$, or
$\fp\fq(2)$, or $\fp\fs\fq(2)$ in
$\fvect(1;\un|1)$ for any $\un$.

\ssec{Prolong of the nonpositive part of $\varphi(\fp\fq(n+1))$ for
$n+1>2$} Let $\id$ denote the tautological $n|n$-dimensional
$\fq(n)$-module. Let us compute the Cartan prolongs of the nonpositive
part of $\varphi(\fp\fq(n+1))$, i.e., compute $(\id, \fq(n)\oplus
\Kee\cdot J)_{*,\un}$ for $p\ne 2$ and $(\id, \fq(n))_{*,\un}$ for $p= 2$.

\sssbegin{Theorem}\label{NewTh} If $p=2$, then $(\id,
\fq(n))_{*,\un}=\fq(\fvect(n;\underline{\One}|0))$.

If $p\neq 2$, then $(\id,
\fq(n))_i=0$ for $i\geq 1$.
\end{Theorem}

\begin{proof} Let $\fh:=(\id, \fq(n))_{*,\un}$ and $X=\mathop{\sum}
(F_i(z,\xi)\pd{z_i}+G_i(z,\xi)\pd{\xi_i})\in \fh$. Let us
single out the image $\fg$
of the 0-th component of $\fh$ in $\fvect(n;\un|n)$  by means of a system of linear differential equations
using operators $\pd{z_i}$ and
$\widetilde{\pd{\xi_i}}=\pd{\xi_i}\cdot\Pty$, where
$\Pty(x):=(-1)^{p(x)}x$ is the parity operator. The system of linear
differential equations obtained defines the (generalized)
Cartan prolong, see \cite{Shch}. The image of $\fq(n)$ is
described by the second line of eqs.~\eqref{.12} or, equivalently,
by means of the following equations:
\begin{equation}\label{deq}
\pd{z_j}F_i=-(-1)^{p(G_i)}\pd{\xi_j}G_i; \quad
\pd{z_j}G_i=-(-1)^{p(F_i)}\pd{\xi_j}F_i \text{ for any } i,j\in
1,\dots, n.
\end{equation}

Since the indices $i$ and $j$ in equations \eqref{deq} are independent of each other, it suffices to find all pairs of
functions $(F,G)$ satisfying the system of equations
\begin{equation}\label{deq1}
\pd{z_j}F=-(-1)^{p(G)}\pd{\xi_j}G; \quad
\pd{z_j}G=-(-1)^{p(F)}\pd{\xi_j}F \text{ for any } j\in 1,\dots, n.
\end{equation}
In these terms, the Lie superalgebra $\fh$ is spanned by vector
fields of the form $F\pd{z_i}+G\pd{\xi_i}$ and
$G\pd{z_i}+F\pd{\xi_i}$ for any $i=1,\dots n$, minding conditions~\eqref{deq1}.

Let the function $F$ be a sum of monomials, one of which is
$Az_i^kf$ (for $p=0$) or $Az_i^{(k)}f$ (for $p>0$), where $f$ does
not depend on $z_i$, and $A\in\Kee$. Because
$\pd{z_i}F=-(-1)^{p(G)}\pd{\xi_i}G$, the monomial $f$ should not
depend on $\xi_i$ either, and the function $G$ should contain a
monomial $z_i^{k-1}\xi_if$ (for $p=0$) or $z_i^{(k-1)}\xi_if$ (for
$p>0$) with, perhaps, different coefficient. But since all equations
are symmetric with respect to the transposition
$F\longleftrightarrow G$, the degree $k-1$ should be equal to $0$.
Thus, the degree of each of the functions $F$ and $G$ relative each
even indeterminate should not exceed $1$.

Let us explicitly describe the component $\fh_1$.

If the vector field $X\in\fh_1$ is homogeneous with respect to
parity, then the functions $F$ and $G$ should be of opposite parity,
and since the equations are symmetric with respect to the
transposition $F\longleftrightarrow G$, let us assume, for
definiteness sake, that \begin{equation}\label{FevGodd} \text{$F$ is even and $G$
is odd.} \end{equation} Therefore, these functions are of the form:
\begin{equation*}\label{g1}
F=\mathop{\sum}\limits_{i<j}(a_{ij}z_iz_j+b_{ij}\xi_i\xi_j), \quad
G=\mathop{\sum}\limits_{i\ne j}c_{ij}z_i\xi_j.
\end{equation*}
Hence
\begin{equation*}\label{F&G1}
\renewcommand{\arraystretch}{1.4}
\begin{array}{ll}
\pd{z_i}
F=\mathop{\sum}\limits_{j>i}a_{ij}z_j+\mathop{\sum}\limits_{j<i}a_{ji}z_j,
&\pd{\xi_i}G=
\mathop{\sum}\limits_{i\ne j}c_{ji}z_j, \\
\pd{\xi_i}F=\mathop{\sum}\limits_{j>i}b_{ij}\xi_j-
\mathop{\sum}\limits_{j<i}b_{ji}\xi_j, &
\pd{z_i}G=\mathop{\sum}\limits_{i\ne j}c_{ij}\xi_j.\\
\end{array}\end{equation*}

Now, eqs.~\eqref{deq1} imply that
\begin{equation*}\label{eqs}
c_{ij}=-b_{ij}=a_{ij} \text{~~and~~} c_{ji}=b_{ij}=a_{ij} \text{~~
for~~} i<j.
\end{equation*}

\underline{If $p\ne 2$}, these conditions contradict each other, and
hence
\begin{equation*}\label{eqsCor}
\text{$(\id, \fq(n))_i=0$ for $i\geq 1$}.
\end{equation*}

\underline{If $p=2$}, eqs.~\eqref{deq1} are free of contradictions,
generate $\binom n2$ pairs of functions $(F,G)$, and, accordingly,
$\binom n2 (n|n)$ vector fields of the form:
\[
\renewcommand{\arraystretch}{1.4}
\begin{array}{l}
X_\ev=(z_iz_j+\xi_i\xi_j)\pd{z_k}+(z_i\xi_j+\xi_iz_j)\pd{\xi_k},
\\
X_\od=(z_iz_j+\xi_i\xi_j)\pd{\xi_k}+(z_i\xi_j+\xi_iz_j)\pd{z_k}.\\
\end{array}\]

Let us describe now the general shape of the pairs of functions
$(F,G)$ satisfying eq.~\eqref{deq1} and defining the vector fields
$X\in\fh$. Let the field $X$ be homogeneous with respect to the
parity and $\Zee$-grading. Then the functions $F$ and $G$ are also
homogeneous with respect to the parity and $\Zee$-grading and their
parities are opposite; recall eq.~\eqref{FevGodd}.

Let $F$ contain, as a summand, a monomial $\alpha\xi_i\xi_j\cdot f$,
where $\alpha\in\Kee$ and $f$ is a monomial independent of $\xi_i$
and $\xi_j$. By differentiating $F$ with respect to $\xi_i$ we see,
by virtue of eqs.~\eqref{deq1}, that the function $G$ contains the
monomial $\alpha z_i\xi_j\cdot f$ as a summand. Differentiating now
the function $G$ with respect to $\xi_j$ we see that $F$ contains
also the term $\alpha z_iz_j\cdot f$.

Consequently applying this
procedure several times we conclude that $F$ must contain monomials
depending on $z$ only, and the sum of all these monomials completely
determines the whole pair $(F,G)$. Let us call this sum \textbf{the
main part} of $F$ and denote it $F^{m}$. For $F^{m}=z_1z_2\dots
z_k$, the pair $(F,G)$ is of the form:
\begin{equation*}\label{F&G}
\renewcommand{\arraystretch}{1.4}
\begin{array}{ll}
F=& z_1z_2\dots z_k+\mathop{\sum}\limits_{1\le i_1<i_2\le k}
\pd{z_{i_1}}\pd{z_{i_2}}(z_1z_2\dots z_k)\xi_{i_1}\xi_{i_2}+\\
& \mathop{\sum}\limits_{1\le i_1<i_2<i_3<i_4\le
k}\pd{z_{i_1}}\pd{z_{i_2}}\pd{z_{i_3}}\pd{z_{i_4}}(z_1z_2\dots
z_k)\xi_{i_1}\xi_{i_2}\xi_{i_3}\xi_{i_4}+\dots.\\

G=& \mathop{\sum}\limits_{i}\pd{z_{i}}(z_1z_2\dots z_k)\xi_{i}+
\mathop{\sum}\limits_{1\le i_1<i_2<i_3\le
k}\pd{z_{i_1}}\pd{z_{i_2}}\pd{z_{i_3}}(z_1z_2\dots
z_k)\xi_{i_1}\xi_{i_2}\xi_{i_3}+\dots.\\
\end{array}\end{equation*}

Let
\begin{equation}\label{mainP}
X=\mathop{\sum}\limits_i
(F_i(z,\xi)\pd{z_i}+G_i(z,\xi)\pd{\xi_i})\longmapsto
X^m=\mathop{\sum}\limits_i F^m_i\pd{z_i}
\end{equation}
be the passage to the
\textbf{main part} $X^m$ of a given even vector $X$.
Clearly, $X^m$ depends only on $z$, and if $X,Y$
are even fields, then
\[
[X^m,Y^m]=[X,Y]^m.
\]
This means that $\fh_\ev\simeq\fvect(n;\underline{\One}|0)$. This isomorphism
is given by the map \eqref{mainP}.

If the field $X$ is odd, then the field $X$, and its main part, are,
due to our assumption \eqref{FevGodd}, of the form
\begin{equation}\label{mainP1}
X=\mathop{\sum}\limits_i (G_i(z,\xi)\pd{z_i}+F_i(z,\xi)\pd{\xi_i})\longmapsto
X^m=\mathop{\sum}\limits_iF^m_i(z)\pd{\xi_i}.
\end{equation}
If now $Y$ is an
even field, then
\[
[Y^m, X^m]=[Y,X]^m,
\]
i.e., as $\fh_\ev$-module, $\fh_\od$ is the adjoint one. It is easy to
verify that $\fh=\fq(\fvect(n;\underline{\One}|0))$.\end{proof}

\sssbegin{Theorem}\label{NL} If $p\ne 2$ and $n>1$, then
\[
\fh:=(\varphi(\fg_{-1}), \varphi(\fg_0))_{*,\un}= (\id,\fq(n)\oplus
\Kee\cdot J)_{*,\un}=\varphi(\fq(n+1)).
\]
\end{Theorem}

\begin{proof} As we repeatedly
demonstrated describing the exceptional simple vectorial Lie
superalgebras, see \cite{Sh5, Sh14}, having added one odd central
element to the 0-th component (this is precisely what happens under
the homomorphism $\varphi$) we can add  to $\fh_1$ either the $\fg_0$-module
$\Pi(\fg_{-1}^*)$, or nothing. Since $(\id, \fq(n))_1=0$
by Theorem \ref{NewTh}, we see that $\fh_1=\varphi(\fg_1)$, where
$\fg_1$ is defined in eq.~\eqref{.11}.

Let us describe the component $\fh_2$. Once again, let us write the
system of differential equations on coordinates of the vector field
$X=\sum(F_i\del_{z_i}+G_i\del_{\xi_i})\in\fh$. Because the
difference from the case of Theorem \ref{NewTh} consists of one
basis element in the 0-th component only, these equations do not
differ much from eqs.~\eqref{deq}. More precisely, the first group
of eqs.~\eqref{deq} remains the same; the second group is the same
for all $i=1,\dots, n$ and all $j\ne i$, whereas $n$ equations of
the second groups for $j=i$ turn into the following $n-1$ equations:
\[
\del_{z_i}G_i+(-1)^{p(F_i)}\del_{\xi_i}F_i=
\del_{z_1}G_1+(-1)^{p(F_1)}\del_{\xi_1}F_1, \text{~~where
$i=2,\dots, n$}.
\]

Let now $X\in \fh_2$ be even. Then the degrees of all (even)
functions $F_i$ and all (odd) functions $G_i$ are equal to 3, i.e.,
\[
\renewcommand{\arraystretch}{1.4}
\begin{array}{ll}F_i=f_i(z)+\sum_{j,k,l} c_i^{jkl}z_j\xi_k\xi_l,
\text{~~where $c_i^{jkl}\in\Kee$ and $\deg(f_i(z))=3$}, \\
G_i=\sum_jg_i^j(z)\xi_j+\sum_{j,k,l} d_i^{jkl}\xi_j\xi_k\xi_l,
\text{~~where $d_i^{jkl}\in\Kee$ and $\deg(g_i^j(z))=2$}.
\end{array}
\]

The following are all equations for the $F_i$ and $G_i$, their parities being
taken into account:
\begin{equation}\label{FG1}
\begin{array}{ll}
\del_{z_j}F_i=\del_{\xi_j}G_i&\text{for any~~} i,j=1,\dots
n\\

\del_{z_j}G_i=-\del_{\xi_j}F_i &\text{for any } i=1,\dots n, j\ne
i;\\

\del_{z_i}G_i+\del_{\xi_j}F_i=\del_{z_1}G_1+\del_{\xi_1}F_1
&\text{for any } i=2,\dots n.
\end{array}
\end{equation}
The first two lines of eqs.~\eqref{FG1} imply:
\[
\del_{\xi_j}\del_{\xi_k}(G_i)=\del_{\xi_j}\del_{z_k}(F_i)=
\del_{z_k}\del_{\xi_j}(F_i)=-\del_{z_k}\del_{z_j}(G_i)\text{~~ for
$j,k\ne i$}.
\]
In particular, for $j=k\ne i$, we get
\[
\frac{\del^2 G_i}{\del z_j^2 }=0 \text{ for all } j\ne i.
\]
Analogously, $\del_{\xi_k}\del_{\xi_j}(G_i)=
-\del_{z_j}\del_{z_k}(G_i)$, and because partial derivatives with
respect to the even variables commute, whereas those with respect to
the odd ones anticommute, we get:
\begin{equation}\label{G4}
\frac{\del^2 G_i}{\del z_j \del z_k}=\frac{\del^2 G_i}{\del \xi_j
\del \xi_k} =0 \text{~~for all } j,k\ne i,
\end{equation}
implying that $d_i^{jkl}=0$, and all functions $g_i^j$ are linear
combinations of monomials $z_i^{(2)}$ and $z_iz_k$, where $k=1,\dots
n$, and $k\ne i$.

We similarly deduce that for $j,k\ne i$, we have
\begin{equation}\label{F5}
\frac{\del^2 F_i}{\del z_j^2 }=\frac{\del^2 F_i}{\del z_j \del
z_k}=\frac{\del^2 F_i}{\del \xi_j \del \xi_k} =0,
\end{equation}
implying that among the coefficients $c_i^{jkl}$ only those of the
form $c_i^{jil}$ can be nonzero, and the function $f_i$ is a linear
combination of the monomials $z_i^{(3)}$ and $z_i^{(2)}z_j$.

Finally, eqs.~\eqref{F5} for $j,k\ne i$ imply that
\[
\frac{\del^2 G_i}{\del z_j \del \xi_k}= \frac{\del^2 F_i}{\del z_j
\del z_k}=0.
\]
Therefore the function $G_i$ is of the form
\begin{equation*}\label{G6}
G_i=z_i^{(2)}\sum_{`\leq j\leq n} a^j_i\xi_j + z_i\xi_i\sum_{j\ne
i}b^j_iz_j.
\end{equation*}

Similarly, eq.~\eqref{G4} for $j,k\ne i$ implies that
\[
\frac{\del^2 F_i}{\del \xi_j \del z_k}= \frac{\del^2 G_i}{\del \xi_j
\del \xi_k}=0\text{~~~~and~~~~}
\frac{\del^2 F_i}{\del \xi_i \del z_i}= \frac{\del^2 G_i}{\del
\xi_i^2}=0.
\]
These equalities imply that $c_i^{jkl}=0$ for any $i,j,k$, and the
function $F_i$ is of the form
\[
F_i=\alpha_iz_i^{(3)}+z_i^{(2)}\sum_{j\ne i}\beta^j_iz_j.
\]
Now, from eq.~\eqref{FG1} for $j\ne i$ we deduce that $a_i^j=\beta_i^j$,
whereas for $j=i$ we have $\beta_i^k=b_i^k$ and $\alpha_i=a_i^i$.
Because $\del_{\xi_j} F_i=0$ for all $j$, the equality \eqref{FG1}
implies $b_j^i=0$ for all $i,j$. As a result,
$F_i=\alpha_iz_i^{(3)}$, $G_i=\alpha_iz_i^{(2)}\xi_i$ and
\[
\frac{\del F_i}{\del_{\xi_i}}+\frac{\del
G_i}{\del_{z_i}}=\alpha_iz_i\xi_i.
\]
But since this sum should not, thanks to the third line of eqs.~\eqref{FG1}, depend on $i$, we conclude that $\alpha_i=0$ for any
$i$, and hence $X=0$.

For any odd vector field $X$ we similarly get
$\fh_i=0$ for any $i\geq 2$. \end{proof}

\ssec{Queerification and Cartan prolongation are commuting
operations for  $p=2$} Let $\fg_0$ be a Lie algebra with a $2$-structure,
$\fg_{-1}$ a restricted $\fg_0$-module, and
$\fg=(\fg_{-1},\fg_0)_{*,\underline{\One}}$ the Cartan prolong.
Now, consider:
\[
\renewcommand{\arraystretch}{1.4}
\begin{tabular}{l}
$\fq\fg=\fg\oplus \Pi(\fg)$, the queerification of the prolong,\\

$\fq\fg_0=\fg_0\oplus\Pi(\fg_0)$, the queerification of the $0$-th
component of the prolong,\\

$\fq\fg_{-1}=\fg_{-1}\oplus\Pi(\fg_{-1})$, the queerification of the
$-1$-st
component of the prolong, \\

$\fg\fq:=(\fq\fg_0, \fq\fg_{-1})_{*,\un}=\oplus \fg\fq_{k,\un}$, the prolong  
of
the queerifications.\\
\end{tabular}
\]
\sssbegin{Theorem}\label{NewTh2} The Lie superalgebras $\fq\fg$ and
$\fg\fq$ coincide: $\fq\fg=\fg\fq$.
\end{Theorem}

\begin{proof} First of all, recall, see \cite{BGLLS1}, that the coordinates of the shearing
vector that can not exceed a certain value are said to be
\textit{critical} and notice that if $\dim \fg_{-1}=n$, then $\fq\fg_0\subset \fq(n)$.
Now, due to Theorem~\ref{NewTh},
$\fg\fq\subset(\id,\fq(n))_{*,\un}=\fq(\fvect(n;\underline{\One}|0))$, and hence all
coordinates of the shearing vector for $\fg\fq$ are critical: $\un=\underline{\One}$.

Observe that the nonpositive parts of the $\Zee$-graded Lie
superalgebra $\fq\fg$ and $\fg\fq$ coincide by their definitions.
Because the Cartan prolongation with a fixed $\un$, again by
definition, is the maximal transitive Lie (super)algebra with
given nonpositive part and $\un$, we get an embedding
$\fq\fg\subset \fg\fq$.

On the other hand, the map $X\mapsto X^m$ constructed in eqs.~
\eqref{mainP}, \eqref{mainP1} determines the embedding
$\fg\fq\tto\fq(\fvect(n;\underline{\One}|0))$. If $X^m$ is even, then $X^m\in
\fg$ by definition of the Cartan prolong and equality of the
shearing vectors for $\fg$ and $\fg\fq$. If $X^m$ is odd, then
$X^m\in \Pi(\fg)$ due to the commutation relations in
$\fq(\fvect(n;\underline{\One}|0))$. Hence, $\fg\fq$ is embedded into $\fq\fg$
and $\fq\fg=\fg\fq$.
\end{proof}

\section{Examples of simple Lie superalgebras obtained by ``Method 2"}\label{Z2gra}

\ssec{How would ``method 2" work if it were defined for $p\neq2$}\label{ssMet2} Over any ground field, ``method 2" applied to $\fsl(n)$ boils down to the following: we declare several pairs of (positive and the
corresponding to them negative) Chevalley generators odd and simultaneously change the corresponding
diagonal elements of the Cartan matrix (replace 2 with a 0). Having changed the Cartan matrix and parities of generators we also change the
defining relations. The Lie superalgebras $\fsl(a|b)$, where $a+b=n$, thus obtained are simple as long
as $a\neq b$, whereas $\fsl(a|a)$ is simple modulo center.

Over $\Cee$ and $\Kee$ with $\Char \Kee>2$, for any simple Lie algebra $\fg\not\simeq\fsl(n)$, having similarly
declared any pair of Chevalley generators odd and factorized by the
ideal of relations that replace Serre relations (for their description, see \cite{BGLL} and references therein) we get a simple (perhaps, modulo
center) Lie superalgebra of infinite dimension; moreover, it is of
rather fast growth as a $\Zee$-graded algebra, cf. \cite{CCLL, GL3}. Under
such superization not only relations that replace Serre ones (the ones
between root vectors of the same sign) change their form and
affect the dimension of the quotient, but also
Cartan matrix and ensuing weight relations become modified: compair
$\fsl(2)$ whose Cartan matrix is $(2)$ and weight relations are $[H, X_\pm]=\pm2X_\pm$
with $\fsl(1|1)$ whose Cartan matrix is $(0)$ and weight relations are $[H, X_\pm]=0$.

\ssec{``Method 2" for vectorial Lie algebras: $p=2$ and $\un=\One$}\label{ssMet2One} In \cite{LeP}, there are given examples of superization of simple
vectorial Lie algebras by what we call here ``method 2" for the
cases where several of coordinates of $\un$ are equal to 1. In this case, one can
consider the corresponding indeterminates odd. In this way we get
$\fh_\Pi(a;\tilde \un|b)$ from $\fh_\Pi(a+b;\un)$ (resp.
$\fle(a;\tilde \un)$ from $\fh_\Pi(2a;\un)$), where $\tilde \un$ is
the part of $\un$ corresponding to the even indeterminates. The same
concerns other series and exceptional examples, see \cite{BGLLS1,
BGLLS2}. These superizations are expected, to an extent. Numerous
examples of simple Lie superalgebras given below are totally new
(except for ``occasional isomorphisms" that might occur for a small
number of indeterminates): these Lie superalgebras are \textit{$\Zee$-graded
and this grading is compatible with $\Zee/2$-grading} whereas over
$\Cee$ only the following simple Lie superalgebras have italicized property:
the series $\fk(1|m)$, and exceptions $\fk\fas$, $\fv\fle(4|3;
K)=\fv\fle(3|6)$, $\fk\fsle(9|6; K)=\fk\fs\fle(5|10)$ and
$\fm\fb(4|5; K)=\fm\fb(3|8)$, see \cite{LSh, Sh14}.

In the cases below, we have to add squares of elements of degree
$-1$ (and $-3$, if there are any; the only simple
infinite-dimensional vectorial Lie algebras with Weisfeiler gradings
known to us for $p=2$ are of depth at most 3) but we do not have to add squares of elements positive
degrees, see Proposition~ \ref{PrNonPos}.

``Method 2", first consciously used in
\cite{BGLLS}, being applied to the deforms of ``standard" vectorial Lie algebras,
in particular, Kaplansky, Eick and Skryabin algebras, see
\cite{BGLLS,Ei,SkT1} yields many new simple Lie superalgebras.

\ssec{``Method 2" applied to simple $\Zee/2$-graded vectorial Lie algebras for which
all coordinates of $\un$ corresponding to \lq\lq odd\rq\rq\ indeterminates are
$>1$} We introduce {$\Zee/2$-gradings} as the $\Zee$-gradings modulo
2. Observe that the desuperized finite-dimensional simple vectorial Lie superalgebra
$\textbf{F}(\fG)(\un)$ over $\Kee$ may have more (types of)
Weisfeiler $\Zee$-gradings than its infinite-dimensional namesake $\fG$ over $\Cee$
has.

\sssec{The simple serial vectorial Lie algebras $\fg=\fg(\un)$ in the standard
$\Zee$-grading $\st$} Recall that in the standard $\Zee$-grading
$\st$ the degree of each indeterminate is equal to 1, except for the
contact case, where the degree of the ``time" indeterminate is equal to 2.

In certain cases, we have to take not the Cartan prolong listed in
the titles of subsections below but its derived since we superize
the \textbf{simple} Lie algebras; but since the simple Lie
superalgebra $\fs(\fg, \gr)$, where $\fg=\fg(\un)$, see Statement~\ref{Statement}, is completely
determined by $\un$, and the elements of degree $\leq 0$ which for $\fg^{(i)}$ coincide
with those of $\fg$, see Proposition~\ref{PrNonPos}, we do not indicate such subtleties.


\paragraph{$\fvect^{(1)}(1;\un)$} We have to add the square of the basis
element of degree $-1$ (obviously, $N>1$) and the square of $x^{(2^{N-1})}\del$; the Lie superalgebra
$\fs(\fvect^{(1)}(1;\un), \st)$ is, clearly, isomorphic to
$\fk(1;\underline{N-1}|1)$. (For far from obvious description of $\fk(1;\un|1)$, see \cite{BGLLS2}.)

\paragraph{$\fvect(n;N)$ for $n>1$, and $\fsvect(n;N)$ for $n>2$;
various Hamiltonian series (all in the standard grading $\st$), and filtered deforms thereof} Set
\be\label{squares}
\fs(\fg, \st)_{-2}:=\fg_{-1}^{[2]}:=\Span(\del_i^2\mid N_i>1).
\ee
In particular, if $\un=\One$, then $\fs(\fg, \st)_{-2}=0$. 

\paragraph{$\fk(2n+1;N)$ for $n>0$ in the standard grading $\st$, and filtered deforms of $\fk(2n+1;N)$}
Having added the trivial $\fg_0$-module $\fg_{-1}^{[2]}:=\Span(\del_i^2\mid N_i>1)$ we get
$\fs(\fg, \st)_{-2}=\fg_{-1}^{[2]} \oplus \fg_{-2}$.

\sssec{The simple vectorial Lie algebras in nonstandard
$\Zee$-gradings} In addition to the ``standard" gradings there are
finitely many (classes of) \textit{nonstandard} $\Zee$-gradings
associated with maximal subalgebras (of finite codimension if the
algebra is infinite-dimensio\-nal) and the corresponding
\textit{Weisfeiler} filtration. Generally, these gradings are given
by setting $\deg x_i=0$ for $i=0, 1, \dots n-1$ of the $n$
indeterminates; for the classification of nonstandard gradings of simple vectorial Lie superalgebras
 over $\Cee$, see \cite{LSh}.

\paragraph{$\fii\fr(9;\un)$ and \textbf{F}($\fv\fle(3;\un|6))$ with a grading $\gr$ described in \cite{BGLLS1}}
Having added the trivial $\fg_0$-module $\fg_{-1}^{[2]}:=\Span(\del_i^2\mid N_i>1)$ we get
$\fs(\fg, \gr)_{-2}=\fg_{-1}^{[2]}\oplus \fg_{-2}$.

\paragraph{\textbf{F}($\fk\fs\fle(5;\un|10))$ and \textbf{F}($\fm\fb(3;\un|8))$
described in \cite{BGLLS2}} These Lie algebras $\fg$ have $-3$rd
components, so their superizations corresponding to such gradings, denote them $\gr3$,  are of depth 6, the squares of the basis elements of degree $-3$ \textbf{corresponding to the indeterminates whose heights are $>1$} span the trivial
$\fg_0$-module we denote hereafter (for other algebras of depth 3 as well)
by $\fg_{-3}^{[2]}$. It is easy to see
that if $\bar\fg $ is the restricted closure of $\fg$, then
$[\bar\fg , \bar\fg]\subset[\fg,\fg]$. In particular, $[\bar\fg,
\fg]\subset\fg$. Therefore, $\fs(\fg,\gr3)_{-4}=\fs(\fg,\gr3)_{-5}=0$
since $\fg$ is of depth 3.

\ssbegin{Proposition}\label{PrNonPos} Let $\fg=\fg(\un)$ be the generalized Cartan prolong. To obtain the Lie superalgebra $\fs(\fg, \gr)$, we only have to add to $\fg$ the squares of basis elements of negative degrees \textbf{corresponding to the indeterminates whose heights are $>1$}.
\end{Proposition}

Recall that in all known examples of simple vectorial Lie superalgebras except $\fv\fle$, for
which one of the coordinates of the shearing vector $\un$ can not exceed
2, if $\un_i$ can be $>1$ it can be however big, see \cite{BGLLS2}.  The coordinates of the
shearing vector with ``in-built" bounds are called \textit{critical}.

\begin{proof} Let the depth of $\fg$ be equal to 1. Let $\fg=(\fg_{-1}, \fg_0)_{*,\un}$
be the infinite-dimensional Cartan prolong corresponding to the
shearing vector $\un$ without any constraints imposed on its
coordinates, so all noncritical coordinates of $\un$ are equal to
$\infty$. Then $\fg$ is the maximal transitive $\Zee$-graded Lie
algebra with the given nonpositive part. In particular, the
maximality implies that the positive part is a restricted Lie
algebra. So, to cook $\fs(\fg, gr)$, it suffices to add squares of the
basis elements of degree $-1$ to $\fg$.

Take the new nonpositive part $\fG:=\mathop{\oplus}\limits_{-2\leq
k\leq 0}\fG_k$ with $\fG_{-1}:=\fg_{-1}$ considered odd,
$\fG_{0}:=\fg_{0}$, and $\fG_{-2}:=\fg_{-1}^{[2]}$, see eq.~\eqref{squares}, added; let
$\fG:=\mathop{\oplus}\limits_{-2\leq k}\fG_k$ be the Cartan prolong
of the new nonpositive part. If we desuperize $\fG$, then in the Lie
algebra $\textbf{F}(\fG)$ obtained, the $-1$st component is
commutative, and $\textbf{F}(\fG)_{-2}$ is a trivial
$\textbf{F}(\fG)_{0}=\fg_0$-module. This means that
$\fii:=\mathop{\oplus}\limits_{-1\leq k}\fG_k$ is an ideal in
$\textbf{F}(\fG)$ with the nonpositive part equal to $\fg_{-1}\oplus
\fg_{0}$. But $\fg$ is the maximal Lie algebra with such nonpositive
part.

The arguments for depth $>1$ are similar.\end{proof}

\sssec{Example} Apply method 2 to
$\fg=\fvect(n;\un)=\Span(f_i(u)\del_{u_i})_{i=1}^n$ in its standard $\Zee$-grading, where
$u=(u_1,\dots, u_n)$. We get
\[
\fG_{-2}=\Span(\del_{x_i})_{i=1}^n,\ \ \ \ 
\fG_{-1}=\Span(D_i=\del_{\xi_i}+\xi_i\del_{x_i})_{i=1}^n,
\]
and
hence, as spaces, $\fG=\Span(f_i(x,\xi)D_i)_{i=1}^n\oplus\Span(\del_j)_{j=1}^n$. The one-to-one
correspondence between $\fg$ and $\fii$, and the description of $f(u)$ as $f(x,\xi)$, are as follows:
\[
\begin{array}{l}
\del_{u_i} \longleftrightarrow D_i, \quad u_i^{(2k)}
\longleftrightarrow x_i^{(k)},\quad u_i^{(2k+1)} \longleftrightarrow
x_i^{(k)}\xi_i. \end{array}
\]

\section{Deform with an odd parameter and its
desuperization}\label{SoddDef}

\ssec{Simple Lie superalgebras}\label{simpleLie} In the functorial
approach, a superspace $\fg$ is a
\textit{Lie superalgebra} if $\underline{\fg}(C):=(\fg\otimes
C)_\ev$ is a Lie algebra for any supercommutative superalgebra $C$
and any superalgebra homomorphism $C_1\tto C_2$ induces a Lie
algebra homomorphism $\underline{\fg}(C_1)\tto\underline{\fg}(C_2)$
with the natural properties of the through and identity maps.

A given Lie superalgebra
$\fg$ is said to be \textit{simple} if $\dim\fg>1$ and it has no
proper ideals. The ideal is also defined functorially. 

\ssec{Deformations of Lie superalgebras}\label{DefLie}
Infinitesimally, the deformation with parameter $t$ of the
multiplication in a given Lie superalgebra $\fg$ is defined by the
expression \eqref{defOfdef}:
\begin{equation}\label{defOfdef}
\begin{cases}
{}[x,y]+tc(x,y)&\text{for $x$ not proportional to $y$ if both are odd},\\
x^2+tc(x,x)&\text{for $x=y$ odd},
\end{cases}
\end{equation}
where $c$ is a 2-cocycle with adjoint coefficients and $p(t)=p(c)$,
and $x,y\in\fg$. 

Thus, eq.~\eqref{defOfdef} means that the global deformation
\textit{linearly} depending on the parameter is a transition
\begin{equation*}\label{defParam}
\fg\tto \fg\otimes\begin{cases}
\cO(1;\underline{\One})=\Kee[t;\underline{\One}]\simeq\Kee[t]/(t^2)
&\text{if $p(t)=\ev$}\\
\Lambda(1)=\Lambda[t]&\text{if $p(t)=\od$}.
\end{cases}
\end{equation*}
The desuperization functor \textbf{F} sends $\Lambda[t]$ to
$\cO(1;\One)$. (These spaces only differ by the parity of $t$; this difference can be seen at
the level of (co)homology, more precisely (co)chains, of degree $>1$.)
Evaluating $\cO(1;\One)\tto\Kee$ by sending $t\longmapsto a\in\Kee$
for some $a$ we turn the Lie algebra over $\cO(1;\One)$
into a Lie algebra over $\Kee$, equivalently $t$ can be considered an even parameter.
For examples, see \cite{BGL3}.


\section{Related open problems}

At the conference \url{http://reims.math.cnrs.fr/pevzner/records.html}  in ho\-nor of A.~Ki\-ril\-lov, P.~Etingof and Yu.~Neretin posed questions some of which we interpret as follows:

1) Describe (compare with \cite{Prem}) the algebraic
groups corresponding to the Brown algebra $\fbr(3)$, the
Weisfeiler-Kac algebras $\fwk^{(1)}(3;a)/\fc$
and $\fwk(4;a)$, the orthogonal Lie algebras $\fo^{(i)}_I(2n)/\fc$;
simple vectorial Lie algebras (especially, for $p=2$), see
\cite{BGLLS1,BGLLS2,Ei}, and to deforms thereof, see \cite{WK,BGL2,
BLW,BGL3}.

2) Describe the algebraic supergroups corresponding to superizations
of the Lie algebras of item 1) and those introduced in this paper (combine \cite{Prem} with \cite{FiGa}).
Describe the automorphism groups for the simple Lie algebras, not
considered in \cite{FG} and \cite{W}, and establish isomorphisms
between the deforms with the help of these groups \`a la \cite{KCh,
Ch}. Describe the automorphism supergroups of the simple Lie
superalgebras obtained by the two methods given in
\S\ref{Sproof}.


\end{document}